\documentclass{amsart}
\usepackage{amssymb,amsmath,amsthm,latexsym}
\usepackage{graphics}
\usepackage{amscd}
\usepackage{here}
\theoremstyle{plain}

\parskip=\bigskipamount

\let\egthree=\phi
\let\phi=\varphi
\let\varphi=\egthree
\let\cal=\mathcal

\begin{document}

\title{Word hyperbolic extensions
of surface groups}
\author{Ursula Hamenst\"adt}
\thanks{Partially supported by Sonderforschungsbereich 611}
\date{May 20, 2005}

\begin{abstract}
Let $S$ be a closed surface of genus $g\geq 2$. A finitely
generated group $\Gamma_S$ is an extension of
the fundamental group $\pi_1(S)$ of $S$ if
$\pi_1(S)$ is a normal subgroup of $\Gamma_S$.
We show that the group
$\Gamma_S$ is hyperbolic if and only
if the orbit map for the action of the
quotient group $\Gamma=\Gamma_S/\pi_1(S)$
on the complex of curves is a quasi-isometric embedding.
\end{abstract}
\maketitle

\section{Introduction}

Let $\Gamma$ be a finitely generated group. A finite symmetric set
${\cal G}$ of generators induces a \emph{word norm} $\Vert\,
\Vert$ on $\Gamma$ by defining $\Vert \phi\Vert$ to be the
smallest length of a word in the generating set ${\cal G}$ which
represents $\phi$. For $\phi,\psi\in \Gamma$ let $d(\phi,\psi)=
\Vert \phi^{-1}\psi\Vert$; then $d$ is a distance function on
$\Gamma$ which is invariant under the left action of $\Gamma$ on
itself. Any two such distance functions on $\Gamma$ are
bilipschitz equivalent. The group $\Gamma$ is called \emph{word
hyperbolic} if equipped with the distance induced by one (and
hence every) word norm, $\Gamma$ is a hyperbolic metric space.

In this note we are interested in
word hyperbolic groups which
are \emph{extensions} of
the fundamental group
$\pi_1(S)$ of a closed orientable surface $S$
of genus $g\geq 2$. By definition,
this means that such a group $\Gamma_S$ contains
$\pi_1(S)$
as a normal subgroup.
Our main goal is to give a geometric
characterization of such groups via
the action of the quotient group
$\Gamma=\Gamma_S/\pi_1(S)$
on the \emph{complex of curves} for
the surface $S$.

Our approach builds on earlier work of
Mosher \cite{Mo96,Mo03} and Farb and Mosher \cite{FM}.
First, recall that
by a classical result of Dehn-Nielsen-Baer (see \cite{I}),
the \emph{extended mapping class group} ${\cal M}_g^0$ of
all isotopy classes of diffeomorphisms of $S$
is just the group of
outer automorphisms of the fundamental
group $\pi_1(S)$ of $S$.
Since the center of $\pi_1(S)$ is trivial we can
identify $\pi_1(S)$ with its group of inner automorphisms
and therefore we obtain an exact sequence
\[1\rightarrow \pi_1(S)\rightarrow
{\rm Aut}(\pi_1(S))\xrightarrow{\Pi} {\cal M}^0_g\rightarrow 1.\]
In particular, for every subgroup $\Gamma$ of ${\cal M}^0_g$ the
pre-image $\Pi^{-1}(\Gamma)$ of $\Gamma$ under the projection
$\Pi$ is an extension of $\pi_1(S)$ with quotient group $\Gamma$.
Vice versa, if $\Gamma_S$ is any group which contains $\pi_1(S)$
as a normal subgroup then the quotient group
$\Gamma=\Gamma_S/\pi_1(S)$ acts as a group of outer automorphisms
on $\pi_1(S)$ and therefore there is a natural homomorphism
$\rho:\Gamma\to {\cal M}^0_g$.

Now consider any finitely generated extension $\Gamma_S$ of
$\pi_1(S)$ with quotient group $\Gamma=\Gamma_S/\pi_1(S)$. If
$\Gamma_S$ is word hyperbolic, then the kernel $K$ of the natural
homomorphism $\rho:\Gamma\to {\cal M}^0_g$ is finite. Namely,
$\Gamma_S$ contains the direct product $\pi_1(S)\times K$ as a
subgroup and no word hyperbolic group can contain the direct
product of two infinite subgroups (see \cite{FM}). As a
consequence, the extension of $\pi_1(S)$ defined by the subgroup
$\rho(\Gamma)$ of ${\cal M}_g^0$ is the quotient of $\Gamma_S$ by
a finite normal subgroup. Since
passing to the quotient by a finite normal subgroup preserves
hyperbolicity, we may assume
without loss of generality that $\Gamma=\Gamma_S/\pi_1(S)$ is a
subgroup of ${\cal M}^0_g$.

The vertices of the \emph{complex of curves} ${\cal C}(S)$
for $S$ are nontrivial free homotopy classes of
simple closed curves on $S$.
The simplices in ${\mathcal C}(S)$
are spanned by collections of such curves which can be realized
disjointly.
In the sequel we restrict our attention to
the one-skeleton of ${\mathcal C}(S)$ which
we denote again by ${\mathcal C}(S)$ by abuse
of notation. Since $g\geq 2$ by assumption,
${\mathcal C}(S)$ is a nontrivial graph
which moreover is connected \cite{Ha}.
However, this graph is locally infinite. Namely,
for every simple closed curve $\alpha$ on $S$ the
surface $S-\alpha$
which we obtain by cutting $S$ open along $\alpha$
contains at least one connected component of Euler characteristic
at most $-2$, and such a component contains infinitely many
pairwise
distinct free homotopy classes of simple closed curves which
viewed as curves in $S$ are disjoint from $\alpha$.

Providing each edge in ${\mathcal C}(S)$ with the standard
euclidean metric of diameter 1 equips the complex of curves
with the structure of a geodesic metric space.
Since ${\mathcal C}(S)$ is not locally finite,
this metric space $({\mathcal C}(S),d)$
is not locally compact. Masur
and Minsky \cite{MM1} showed that nevertheless its geometry
can be understood quite explicitly. Namely, ${\mathcal C}(S)$
is hyperbolic of infinite diameter (see also \cite{B1,H05} for
alternative shorter proofs).
The extended mapping class group
${\cal M}^0_{g}$ of all isotopy classes
of diffeomorphisms of $S$ acts
naturally on ${\cal C}(S)$ as a group of simplicial isometries. In
fact, Ivanov showed that
if $g\not=2$
then ${\cal M}^0_{g}$ is precisely the isometry
group of ${\cal C}(S)$ (see \cite{I} for a sketch of a proof and
for references).

A map $\Phi$ of a finitely generated group
$\Gamma$ into a metric space $(Y,d)$
is called a \emph{quasi-isometric embedding} if for some
(and hence every) choice of a word norm
$\Vert\,\Vert$ for $\Gamma$ there exists a number
$L>1$ such that
\[d(\Phi \psi,\Phi \eta)/L-L\leq \Vert \psi^{-1}\eta\Vert
\leq L d(\Phi \psi,\Phi\eta)+L.\] Note that a quasi-isometric
embedding need not be injective. The following definition extends
the well known notion of a convex cocompact group of isometries of
a simply connected Riemannian manifold of bounded negative
curvature to subgroups of the extended mapping class group, viewed
as the isometry group of the complex of curves.

\bigskip

{\bf Definition:} A finitely generated subgroup
$\Gamma$ of ${\cal M}^0_{g}$ is called \emph{convex
cocompact} if for some $\alpha\in {\cal C}(S)$ the orbit
map $\phi\in \Gamma\to \phi\alpha\in {\cal C}(S)$ is
a quasi-isometric embedding.

\bigskip

For every subgroup $\Gamma$ of ${\cal M}_{g}^0$, the intersection
of $\Gamma$ with the \emph{mapping class group} ${\cal M}_{g}$ of
all isotopy classes of orientation preserving diffeomorphisms is a
subgroup of $\Gamma$ of index at most 2 and hence this group is
quasi-isometric to $\Gamma$. Thus in the sequel we may restrict
our attention to subgroups of ${\cal M}_{g}$.
Since ${\cal C}(S)$
is a hyperbolic geodesic metric space, every convex cocompact
subgroup $\Gamma$ of ${\cal M}_{g}$ is word hyperbolic.

Farb and Mosher \cite{FM} introduce another notion of a convex
cocompact subgroup $\Gamma$ of ${\cal M}_{g}$ via its action on
the \emph{Teichm\"uller space} ${\cal T}_{g}$ of all marked
hyperbolic metrics on $S$. Namely, they define a group
$\Gamma<{\cal M}_{g}$ to be convex cocompact if a $\Gamma$-orbit
in ${\cal T}_{g}$ is \emph{quasi-convex} with respect to the
\emph{Teichm\"uller metric}; this means that for every fixed $h\in
{\cal T}_{g}$ and any two elements $\phi,\psi\in \Gamma$ the
unique Teichm\"uller geodesic connecting $\phi h$ to $\psi h$ is
contained in a uniformly bounded neighborhood of the orbit $\Gamma
h$. Answering among other things a question raised by Farb and
Mosher, we show.

\bigskip

{\bf Theorem:} {\it For a finitely
generated subgroup $\Gamma$
of ${\cal M}_{g}$, the following
are equivalent.
\begin{enumerate}
\item $\Gamma$ is convex cocompact.
\item Some $\Gamma$-orbit on ${\cal T}_{g}$ is quasi-convex.
\item The natural extension of
$\pi_1(S)$ with quotient group $\Gamma$ is word
hyperbolic.
\end{enumerate}
}

\bigskip

The implication $3)\Longrightarrow 2)$ in our theorem is due to
Farb and Mosher \cite{FM} and was the main motivation for this
work. In the particular case that the subgroup $\Gamma$ of ${\cal
M}_g$ is \emph{free}, the reverse implication $2)\Longrightarrow
3)$ is also shown in \cite{FM}. The equivalence of
1) and 2) was independently and at the same time established
by Kent and Leininger \cite{KL05}, with a different proof.

Examples of convex cocompact subgroups of ${\cal M}_{g}$ are
\emph{Schottky groups} which are defined to be \emph{free} convex
cocompact subgroups of ${\cal M}_{g}$. There is an abundance of
such groups: Since every \emph{pseudo-Anosov} element of ${\cal
M}_{g}$ acts with north-south dynamics on the Gromov boundary of
the complex of curves, the classical ping-pong lemma shows that
for any two non-commuting pseudo-Anosov elements $\phi,\psi\in
{\cal M}_{g}$ there are numbers $k\geq 1,\ell\geq 1$ such that the
subgroup of ${\cal M}_{g}$ generated by $\phi^k,\psi^\ell$ is free
and convex cocompact. However, to our knowledge there are no known
examples of convex cocompact groups which are not virtually free.
On the other hand, there are examples of surface-subgroups of
${\cal M}_g$ with interesting geometric properties
\cite{GDH99,LR05}, but these groups contain elements
which are not pseudo-Anosov. Since
the orbit on ${\cal C}(S)$ of an infinite cyclic subgroup
of ${\cal M}_{g}$ generated by an element which
is not pseudo-Anosov is bounded, these groups
are not convex cocompact.

The organization of this paper is as follows. In Section 2, we
define a map $\Psi:{\cal T}_{g}\to {\cal C}(S)$ which is roughly
equivariant with respect to the action of ${\cal M}_{g}$. We
characterize
quasi-geodesics in Teichm\"uller space which are mapped by
$\Psi$ to quasi-geodesics in the complex of curves and deduce from
this as an immediate corollary
the equivalence of 1) and 2) in our theorem. In
Section 3 we give a geometric description of hyperbolic
fibrations with fibre a tree and base a hyperbolic
geodesic space. This is used in Section 4 to show
the equivalence of 1) and 3) in our theorem.

\section{Quasi-geodesics in Teichm\"uller space which project
to quasi-geodesics in the complex of curves}

In this section we
consider an oriented surface
$S$ of genus $g\geq 0$ with $m\geq 0$ punctures. We require
that $S$ is \emph{non-exceptional}, i.e. that $3g-3+m\geq 2$.
The \emph{Teichm\"uller space} ${\cal T}_{g,m}$ of all marked
isometry classes of complete hyperbolic metrics on $S$ of finite
volume is homeomorphic to $\mathbb{R}^{6g-6+2m}$.
The \emph{mapping class group}
${\cal M}_{g,m}$ of all isotopy classes of orientation
preserving diffeomorphisms of $S$ acts
properly discontinuously
on ${\cal T}_{g,m}$
preserving the \emph{Teichm\"uller metric}.
The Teichm\"uller metric is a complete Finsler metric on
${\cal T}_{g,m}$.

The one-skeleton ${\cal C}(S)$ of the
\emph{complex of curves} for $S$ is defined to
be the metric graph
whose vertices are free homotopy classes of simple closed
\emph{essential} curves, i.e. curves which are not contractible
or homotopic into a puncture, and where two
such vertices are connected by an edge of length 1
if and only if they can
be realized disjointly. Since $S$ is non-exceptional
by assumption, the graph ${\cal C}(S)$ is connected.
Moreover, as a metric space it is hyperbolic in the sense of
Gromov \cite{MM1,B1,H05}.

For every marked hyperbolic
metric $h\in {\cal T}_{g,m}$,
every essential free homotopy
class $\alpha$ on $S$ can be represented by a closed
geodesic which is unique up to parametrization.
This geodesic is simple if the free homotopy
class admits a simple representative. The
\emph{$h$-length} $\ell_h(\alpha)$
of the class is defined to be the length of its geodesic representative;
equivalently, $\ell_h(\alpha)$ equals
the minimum of the $h$-lengths of all closed
curves representing the class $\alpha$.

A \emph{pants decomposition} for $S$ is a collection of
$3g-3+m$ pairwise disjoint simple closed
essential curves on $S$
which decompose
$S$ into $2g-2+m$ \emph{pairs of pants}, i.e. planar
surfaces homeomorphic to a three-holed sphere.
By a classical result of Bers
(see \cite{Bu}), there is a number $\chi >0$
only depending on the topological type
of $S$ such that
for every complete hyperbolic metric $h$ on $S$
of finite volume there is
a pants decomposition for $S$ consisting of
simple closed curves of $h$-length at most $\chi$.
Define a map $\Psi:
{\cal T}_{g,m}\to {\cal C}(S)$ by associating
to a complete hyperbolic metric $h$ on $S$
of finite volume an essential simple closed curve
$\Psi(h)\in {\cal C}(S)$ whose $h$-length is at most $\chi$.
By the collar theorem for hyperbolic surfaces
(see \cite{Bu}), the number of intersection points
between any two simple closed geodesics of length
at most $\chi$ is bounded from above by a universal
constant. On the other hand, the distance
between two curves $\alpha,\beta\in {\cal C}(S)$ is
bounded from above by the minimal number of intersection
points between any representatives of $\alpha,\beta$ plus one
\cite{MM1,B1}. Thus
the diameter in ${\cal C}(S)$ of the set
of all simple closed curves of $h$-length at
most $\chi$ is bounded from above by a universal
constant $R>0$ and
the map $\Psi$ is roughly equivariant
with respect to the action of ${\cal M}_{g,m}$. This
means that for every $\phi\in {\cal M}_{g,m}$
and every $h\in {\cal T}_{g,m}$ we have
$d(\phi(\Psi h),\Psi(\phi h))\leq R$.

Let $J\subset \mathbb{R}$ be a closed
connected subset, i.e. either $J$ is a closed interval
or a closed ray or the whole line.
For some $p>1$, a map $\gamma:J\to {\cal C}(S)$ is
called a \emph{$p$-quasi-geodesic} if for all
$s,t\in J$ we have
\[d(\gamma(s),\gamma(t))/p-p\leq \vert s-t\vert
\leq pd(\gamma(s),\gamma(t))+p.\]
The map $\gamma:J\to {\cal C}(S)$ is called
an \emph{unparametrized $p$-quasi-geodesic} if
there is a closed connected set $I\subset\mathbb{R}$ and
a
homeomorphism $\zeta:I\to J$ such that
$\gamma\circ \zeta:I\to {\cal C}(S)$ is a
$p$-quasi-geodesic.
By a result of Masur and Minsky (Theorem 2.6
and Theorem 2.3 of \cite{MM1}, see also
\cite{H05} for a more explicit statement with proof),
there is a number $p>1$ such that
the image under $\Psi$ of \emph{every} Teichm\"uller
geodesic (i.e. every geodesic in ${\cal T}_{g,m}$ with
respect to the Teichm\"uller metric)
is an unparametrized $p$-quasi-geodesic.
However, in general this curve is not
a quasi-geodesic with
its proper parametrization (see \cite{MM1}).

For $\epsilon >0$ let ${\cal T}_{g,m}^{\epsilon}$
be the collection of all hyperbolic
metrics $h\in {\cal T}_{g,m}$ for which
the length of the shortest closed $h$-geodesic
is at least $\epsilon$. Informally we think of
${\cal T}_{g,m}^{\epsilon}$ as the $\epsilon$-thick part
of Teichm\"uller space. The mapping class group preserves
the set ${\cal T}_{g,m}^\epsilon$ and acts on it
cocompactly. Moreover, every ${\cal M}_{g,m}$-invariant
subset of ${\cal T}_{g,m}$ on which ${\cal M}_{g,m}$ acts
cocompactly is contained in ${\cal T}_{g,m}^\epsilon$ for
some $\epsilon >0$.

Define for
$\epsilon >0$ a \emph{quasi-convex curve} in
${\cal T}_{g,m}^\epsilon$ to be a closed subset of ${\cal T}_{g,m}$
whose \emph{Hausdorff distance} to the image of
a geodesic arc $\zeta:J\to {\cal T}_{g,m}^\epsilon$
is at most $1/\epsilon$. Recall that the Hausdorff
distance between two closed subsets $A,B$ of a metric space
is the infimum of all numbers $r>0$ such that
$A$ is contained in the $r$-neighborhood of $B$ and
$B$ is contained in the $r$-neighborhood of $A$.
The main goal of this section is to show the
following result of independent interest.

\bigskip

{\bf Theorem 2.1:} {\it \begin{enumerate}
\item For every $\nu>1$ there is a constant $\epsilon=\epsilon(\nu)>0$
with the following property. Let $J\subset \mathbb{R}$
be a closed connected set of diameter at least $1/\epsilon$ and let
$\gamma:J\to {\cal T}_{g,m}$ be a $\nu$-quasi-geodesic.
If $\Psi\circ \gamma$
is a $\nu$-quasi-geodesic in ${\cal C}(S)$ then $\gamma(J)$ is a
quasi-convex curve in ${\cal T}_{g,m}^\epsilon$.
\item For every $\epsilon >0$ there is a constant
$\nu(\epsilon)>1$ with the following property.
Let $\gamma:J\to {\cal T}_{g,m}$ be a
$1/\epsilon$-quasi-geodesic in ${\cal T}_{g,m}$ whose
image $\gamma(J)$ is a quasi-convex curve
in ${\cal T}_{g,m}^\epsilon$; then
$\Psi\circ \gamma$ is a $\nu(\epsilon)$-quasi-geodesic
in ${\cal C}(S)$.
\end{enumerate}
}

\bigskip

We begin with establishing the second part of our theorem.
For this we need
the following simple no-retraction lemma
for quasi-geodesics in the hyperbolic geodesic metric
space ${\cal C}(S)$.

\bigskip

{\bf Lemma 2.2:} {\it For $p>1$
there is a constant $c=c(p)>0$
with the following property.
Let $\gamma:J\to {\cal C}(S)$ be any unparametrized
$p$-quasi-geodesic; if $t_1<t_2<t_3\in J$
then $d(\gamma(t_1),\gamma(t_3))\geq
d(\gamma(t_1),\gamma(t_2))+d(\gamma(t_2),\gamma(t_3))-c$.}

{\it Proof:} Let $p>1$;
by the definition of an unparametrized
$p$-quasi-geodesic,
it is enough to show the existence of a number
$c>0$ such that for every (parametrized) $p$-quasi-geodesic
$\gamma:[0,n]\to {\cal C}(S)$ and all $0<t<n$
we have $d(\gamma(0),\gamma(n))\geq
d(\gamma(0),\gamma(t))+d(\gamma(t),\gamma(n))-c$.

Since ${\cal C}(S)$ is hyperbolic, there is a constant
$R>0$ only depending on $p$
such that the Hausdorff distance between
every $p$-quasi-geodesic and every geodesic
connecting the
same endpoints is at most $R$.
Let $\gamma:[0,n]\to {\cal C}(S)$ be any $p$-quasi-geodesic and
let $\zeta:[0,m]\to {\cal C}(S)$
be a geodesic connecting $\gamma(0)$ to $\gamma(n)$.
Then for every $t\in (0,n)$ there is a point
$s\in (0,m)$ such that
$d(\gamma(t),\zeta(s))\leq R$.
Thus we have
$d(\gamma(0),\gamma(t))+d(\gamma(t),\gamma(n))\leq
d(\zeta(0),\zeta(s))+d(\zeta(s),\zeta(m))+2R=
d(\gamma(0),\gamma(n))+2R$ which shows the lemma.
\qed

\bigskip

A \emph{geodesic lamination} for a hyperbolic metric
$h\in {\cal T}_{g,m}$ is a \emph{compact} subset of
$S$ foliated by simple $h$-geodesics \cite{CEG}.
A \emph{measured geodesic lamination}
$\mu$ on $S$ is a geodesic lamination together with
a nontrivial transverse invariant measure. An example of a measured
geodesic lamination on $S$ is a simple
closed curve with the transverse counting measure.
The space ${\cal M\cal L}$ of measured geodesic laminations
on $S$ can be equipped with the weak$^*$-topology, and
with this topology, it is homeomorphic to $\mathbb{R}^{6g-6+2m}-\{0\}$.
There is a natural continuous action of
the multiplicative group $(0,\infty)$ of positive
reals on ${\cal M\cal L}$ by scaling, and the quotient of
${\cal M\cal L}$ under this action is the space
${\cal P\cal M\cal L}$ of \emph{projective
measured laminations} which is homeomorphic to
the sphere $S^{6g-7+2m}$. It
can naturally be identified with the projectivized
tangent space of ${\cal T}_{g,m}$ at $h$.
The space ${\cal P\cal M\cal L}$ also
is the boundary of
a compactification of ${\cal T}_{g,m}$, called the
\emph{Thurston boundary} of Teichm\"uller space.
This is used to show.

\bigskip

{\bf Lemma 2.3:} {\it For every $\epsilon >0$ there
is a number $\nu_0=\nu_0(\epsilon)>0$ with the following property.
Let
$\gamma:J\to {\cal T}_{g,m}^\epsilon$ be
a Teichm\"uller geodesic; then the curve
$\Psi\circ \gamma:J\to {\cal C}(S)$ is a
$\nu_0$-quasi-geodesic.}

{\it Proof:} Let $p>1$ be such that the image under $\Psi$
of every Teichm\"uller geodesic is an unparametrized
$p$-quasi-geodesic in ${\cal C}(S)$; such a number exists by the results
of Masur and Minsky \cite{MM1,H05}.
Let $c=c(p)>0$ be as in Lemma 2.2.

We claim that for every $\epsilon >0$
there is a constant $k_0=k_0(\epsilon)>0$ with the
following property. Let $k\geq k_0$ and let
$\gamma:[0,k]\to
{\cal T}_{g,m}^\epsilon$ be a geodesic arc of length
at least $k_0$; then $d(\Psi(\gamma(0)),\Psi(\gamma(k)))
\geq 2c$.

To see that this is the case, we argue by
contradiction and we assume otherwise.
With the number $c$ as above, by Lemma 2.2 there
is then a number $\epsilon >0$ and for every $k>0$
there is a geodesic arc
$\gamma_k:[0,k]\to {\cal T}_{g,m}^\epsilon$ such that
$d(\Psi\gamma_k(0),\Psi\gamma_k(t))\leq 3c$ for
every $t\in [0,k]$. Let $R>0$ be
an upper bound for the diameter in ${\cal C}(S)$
of the set of all simple closed curves whose
length with respect to some fixed metric
$h\in {\cal T}_{g,m}$ is at most $\chi$.
Since the action of ${\cal M}_{g,m}$ on
${\cal T}_{g,m}^\epsilon$ is isometric and cocompact,
via replacing our constant $3c$ by $3c+2R$
we may assume
that the initial points $\gamma_k(0)$
$(k\geq 1)$ of the geodesic arcs $\gamma_k$
are contained
in a fixed compact subset of ${\cal T}_{g,m}^{\epsilon}$.
Thus by passing to a subsequence we may assume that
the geodesics $\gamma_k$ converge locally uniformly
as $k\to \infty$ to
a geodesic $\gamma:[0,\infty)\to {\cal T}_{g,m}^\epsilon$.
By the definition of the map $\Psi$ and continuity of the
length functions on Teichm\"uller space
we then have $d(\Psi\gamma(s),\Psi\gamma(0))
\leq 3c+4R$ for all $s\geq 0$.

Let $\lambda\in {\cal P\cal M\cal L}$
be the projective measured geodesic
lamination which defines the
direction of $\gamma$ at $\gamma(0)$, viewed
as a point in the projectivized tangent space
of ${\cal T}_{g,m}$ at $\gamma(0)$.
Since
$\gamma$ is
cobounded, i.e. it projects into a compact subset of moduli
space ${\cal T}_{g,m}/{\cal M}_{g,m}$,
by a result of Masur \cite{Ma82a}
the lamination $\lambda$ \emph{fills
up} $S$; this means that every simple closed curve on $S$
intersects $\lambda$ transversely. Moreover, $\gamma(t)$
converges as $t\to\infty$ in the Thurston compactification
of ${\cal T}_{g,m}$ to $\lambda$ \cite{Ma82b}.
By the definition
of the Thurston compactification of ${\cal T}_{g,m}$
(see \cite{FLP}), this implies that
the curves $\Psi(\gamma(t))$, viewed
as projective measured laminations,
converge as $t\to \infty$ in ${\cal P\cal M\cal L}$ to $\lambda$.
As a consequence, the curve
$\Psi\circ \gamma$ is an unparametrized
quasi-geodesic in ${\cal C}(S)$ of \emph{infinite} diameter
(see \cite{K}, \cite{H04}) which
is a contradiction and shows our claim.

Now let $n>0$ and
let $\gamma:[0,k_0n]\to {\cal T}_{g,m}^\epsilon$
be any Teichm\"uller geodesic.
The image under $\Psi$ of every
geodesic in ${\cal T}_{g,m}$ is an unparametrized
$p$-quasi-geodesic; thus by the choice of $c$,
for all
$0\leq s\leq t$ we have $d(\Psi\gamma(t),\Psi\gamma(0))\geq
d(\Psi\gamma(t),\Psi\gamma(s))+d(\Psi\gamma(s),\Psi\gamma(0))-c$.
On the other
hand, from our above consideration and the choice
of $k_0$ we conclude that for every
$u<n$ we have
$d(\Psi\gamma(uk_0),\Psi\gamma((u+1)k_0))\geq 2c$ and therefore
$d(\Psi\gamma((u+1)k_0)),\Psi\gamma(s))\geq
d(\Psi\gamma(uk_0),\Psi\gamma(s))+c$ for all $s\leq uk_0$.
Inductively we deduce
that $d(\Psi\gamma(uk_0),\Psi\gamma(vk_0))\geq c\vert u-v\vert$
for all integers
$u,v\leq n$. The map $\Psi:{\cal T}_{g,m}\to
{\cal C}(S)$ is \emph{coarsely Lipschitz} by which we
mean that there is
a constant $a>0$ such that $d( \Psi h,\Psi h^\prime)\leq
ad(h,h^\prime)+a$ for all $h,h^\prime\in {\cal T}_{g,m}$ and where
$d(h,h^\prime)$ denotes the Teichm\"uller distance
between $h$ and $h^\prime$.
Together with above, it follows
that $\Psi\gamma$ is a $\nu_0$-quasi-geodesic for
a constant $\nu_0>0$ only depending on $\epsilon$
(more precisely, we have $c\vert s-t\vert/k_0-k_0a-a\leq
d(\Psi\gamma(s),\Psi\gamma(t))\leq a\vert s-t\vert +a$
for all $s,t\in [0,k_0n]$). This shows the
lemma.
\qed

\bigskip

The following corollary shows the second
part of Theorem 2.1.

\bigskip

{\bf Corollary 2.4:} {\it For every $\epsilon >0$
there is a number $\nu=\nu(\epsilon) >1$ with
the following property. Let
$\gamma:J\to {\cal T}_{g,m}$ be
a $1/\epsilon$-quasi-geodesic
such that $\gamma(J)$ is a quasi-convex curve
in ${\cal T}_{g,m}^\epsilon$; then $\Psi\circ \gamma:J\to {\cal C}(S)$
is a $\nu$-quasi-geodesic.}

{\it Proof:} Let $\epsilon >0$ and let
$\gamma:J\to {\cal T}_{g,m}$ be a $1/\epsilon$-quasi-geodesic
such that $\gamma(J)$ is a quasi-convex curve in
${\cal T}_{g,m}^\epsilon$.
Then there is a Teichm\"uller geodesic
$\zeta:I\to {\cal T}_{g,m}^{\epsilon}$
and a map
$\rho:J\to I$ such that $d(\zeta(\rho(t)),\gamma(t))\leq 1/\epsilon$
for all $t$. Since by assumption the map $\gamma$
is a $1/\epsilon$-quasi-geodesic in ${\cal T}_{g,m}$
and since $\zeta$ realizes
the distance between any of its points, the map $\rho$
is necessarily a $b$-quasi-isometry for a constant
$b>1$ only depending on $\epsilon$. On the other hand,
the map $\Psi$
is coarsely Lipschitz and
therefore the distances $d(\Psi \gamma(t),\Psi(\zeta\circ\rho(t)))$
are bounded from above by a number only depending on
$\epsilon$.
This implies by Lemma 2.3 that
$\Psi\circ \gamma$ is a $\nu$-quasi-geodesic
for a constant $\nu>0$ only depending on $\epsilon$.
\qed

\bigskip

To show the first part of Theorem 2.1,
we begin again with a simple observation.

\bigskip

{\bf Lemma 2.5:} {\it For every $\nu>1$ there is a number
$\epsilon_0=\epsilon_0(\nu)>0$ with the following properties.
Let
$\gamma:[0,n]\to {\cal T}_{g,m}$ be a $\nu$-quasi-geodesic
whose
projection $\Psi\gamma$
to ${\cal C}(S)$ is a $\nu$-quasi-geodesic.
If $n\geq 1/\epsilon_0$ then
$\gamma[0,n]\subset {\cal T}_{g,m}^{\epsilon_0}$.}

{\it Proof:} Let $n>0,\nu>1$ and let
$\gamma:[0,n]\to {\cal T}_{g,m}$ be a $\nu$-quasi-geodesic
such that $\Psi\circ\gamma$ is a $\nu$-quasi-geodesic
in ${\cal C}(S)$.
Then we have $d(\Psi(\gamma(t)),\Psi(\gamma(s)))\geq
\vert s-t\vert/\nu -\nu$ for all $s,t\in [0,n]$. Let
$R>0$ be an upper bound for
the diameter in ${\cal C}(S)$ of the collection
of all simple closed curves on $S$ whose
length with respect to some metric $h\in {\cal T}_{g,m}$
is at most $\chi$ where
as before, $\chi>0$ is determined by Bers' theorem.
Let $[a,b]\subset [0,n]$
be an interval for which there is a simple closed
curve $\alpha\in {\cal C}(S)$ so that $\ell_{\gamma(t)}(\alpha)
\leq\chi$ for all $t\in [a,b]$;
then we have $d((\Psi(\gamma(a)),\alpha)\leq R,
d(\Psi(\gamma(b)),\alpha)\leq R$ and therefore
$\vert b-a\vert\leq 2\nu R+\nu^2$.

Now by a result of Wolpert (see \cite{IT}), for
all $\alpha\in {\cal C}(S)$ and all $h,h^\prime\in
{\cal T}_{g,m}$ the distance between $h$ and $h^\prime$
is at least $\vert\log \ell_h(\alpha)-\log\ell_{h^\prime}(\alpha)
\vert$. Thus if there is a point $t\in [0,n]$ with
$\log(\ell_{\gamma(t)}(\alpha))< \log(\chi)-2\nu R-2\nu^2$
then the $\gamma(s)$-length of $\alpha$ is smaller than
$\chi$ for every $s\in [0,n]$ with
$\vert s-t\vert \leq 2\nu R+2\nu^2$
and consequently by our above consideration, $\Psi\circ \gamma$
is \emph{not} a $\nu$-quasi-geodesic provided that
$n\geq 4\nu R +4\nu^2$.
\qed

\bigskip

Every Teichm\"uller geodesic line
$\gamma:\mathbb{R}\to {\cal T}_{g,m}$ is
defined by a \emph{quadratic differential} on $S$.
More precisely, for each $t\in \mathbb{R}$ there is a holomorphic
quadratic differential $q_t$ on the
Riemann surface $\gamma(t)$ defining
a singular euclidean metric on $S$ in the conformal
class of $\gamma(t)$ and of area one. The differential
$q_t$ and the corresponding piecewise euclidean metric
are determined by the \emph{horizontal} and the
\emph{vertical} foliation of $q_t$. These foliations
have a common
finite set of singular points and are
equipped with a transverse
invariant measure.
For $s\not= t$, the horizontal
foliation for $q_s$ coincides with the horizontal
foliation for $q_t$, but its transverse measure
is obtained from the transverse measure
for $q_t$ by scaling with the factor
$e^{t-s}$. Similarly,
the vertical foliation of $q_s$ coincides with the
vertical foliation of $q_t$, but its transverse measure is
obtained from the transverse measure for $q_t$ by
scaling with the factor $e^{s-t}$.
We use this description of Teichm\"uller geodesics
together with the arguments
in Section 3.9 of \cite{Mo03} to show the first part of
Theorem 2.1.

\bigskip

{\bf Lemma 2.6:} {\it For every $\nu>1$ there is a constant
$\epsilon=\epsilon(\nu)>0$ with the following property. Let
$J\subset \mathbb{R}$ be a closed connected subset of
diameter at least $1/\epsilon$ and let
$\gamma:J\to {\cal T}_{g,m}$ be a $\nu$-quasi-geodesic such that
$\Psi\circ \gamma:J\to {\cal C}(S)$ is a $\nu$-quasi-geodesic
in ${\cal C}(S)$;
then $\gamma(J)$ is a quasi-convex curve in ${\cal T}_{g,m}^\epsilon$.}

{\it Proof:} For $\nu>1$
define a \emph{$\nu$-Lipschitz curve} in ${\cal T}_{g,m}$
to be a $\nu$-Lipschitz map $\gamma:J\to {\cal T}_{g,m}$ with
respect to the standard metric on $\mathbb{R}$ and
the Teichm\"uller metric on ${\cal T}_{g,m}$.
Since ${\cal T}_{g,m}$ is a smooth manifold and
the Teichm\"uller
metric is a complete Finsler metric, every
$\nu$-quasi-geodesic $\gamma:J\to {\cal T}_{g,m}$
can be replaced by a piecewise
geodesic $\zeta:J\to {\cal T}_{g,m}$
which is a
$2\nu$-Lipschitz curve and which satisfies $d(\gamma(t),
\zeta(t))\leq 2\nu$ for all $t\in J$.
Thus it is
enough to show the statement
of the lemma for $\nu$-Lipschitz curves
$\gamma:J\to {\cal T}_{g,m}$ which are $\nu$-quasi-geodesics and
such that
$\Psi\circ\gamma$ is a $\nu$-quasi-geodesic in ${\cal C}(S)$.
In the sequel we also assume that the
diameter $\vert J\vert$ of the set $J$ is
bigger than $1/\epsilon_0$ where
$\epsilon_0=\epsilon_0(\nu)$ is as in
Lemma 2.5; then $\gamma(J)\subset {\cal T}_{g,m}^{\epsilon_0}$.

Since ${\cal C}(S)$ is hyperbolic and $\Psi\circ \gamma$
is a $\nu$-quasi-geodesic
by assumption, there is a geodesic arc in ${\cal C}(S)$
whose Hausdorff distance to
$\Psi\circ \gamma(J)$
is bounded from above by a universal constant.
As a consequence, if $J$ is one-sided infinite, say
if $[0,\infty)\subset J$, then the points
$\Psi(\gamma(t))$ converge as $t\to \infty$ to
a point in the \emph{Gromov boundary}
$\partial {\cal C}(S)$ of
${\cal C}(S)$ (see [BH] for the definition of the
Gromov boundary of a hyperbolic geodesic metric space).
The Gromov boundary
of ${\cal C}(S)$ can naturally be identified with
the space of \emph{minimal}
geodesic laminations on $S$ which fill up $S$, equipped
with a \emph{coarse Hausdorff topology} (see \cite{H04}).
Here a geodesic lamination is minimal if each of its
half-leaves is dense, and it fills up $S$ if it intersects
every essential simple closed curve on $S$ transversely.

A simple closed curve $\alpha\in {\cal C}(S)$ defines
a projective measured lamination which we denote
by $[\alpha]$. Similarly, for
a measured lamination $\lambda\in {\cal M\cal L}$
we denote by $[\lambda]$ the projective class
of $\lambda$. Following Mosher \cite{Mo03},
we say that the projective measured lamination
$[\alpha]$ defined
by a simple closed curve $\alpha\in {\cal C}(S)$
is \emph{realized} at
some $t\in J$ if the length of $\alpha$ with respect
to the metric $\gamma(t)\in {\cal T}_{g,m}$
is at most $\chi$. Note that the number of projective
measured laminations which are realized at
a given point $t\in J$
is uniformly bounded and that $[\Psi(\gamma(t))]$ is
realized at $\gamma(t)$.
Similarly,
we say that the projectivization $[\lambda]$ of a measured
geodesic lamination $\lambda$ is
realized at an infinite ``endpoint'' of $J$ if
the support of $\lambda$ equals
the corresponding endpoint of the quasi-geodesic
$\Psi\gamma(J)$
in the Gromov boundary $\partial {\cal C}(S)$
of ${\cal C}(S)$, viewed as
a minimal geodesic lamination.
The set of projective measured laminations which are
realized at an infinite endpoint of $J$ is a
nonempty closed
subset of ${\cal P\cal M\cal L}$ (see
\cite{K}, \cite{H04}).
We call a projective measured lamination which is realized at
a (finite or infinite) endpoint of $J$
an \emph{endpoint lamination}.

Now $\Psi\gamma$ is a $\nu$-quasi-geodesic in ${\cal C}(S)$
by assumption and the diameter in ${\cal C}(S)$
of the set of all curves of length at most $\chi$
with respect to some fixed hyperbolic metric $h\in {\cal T}_{g,m}$
is bounded from above by a universal constant.
Since any two curves $\alpha,\beta\in {\cal C}(S)$ with
$d(\alpha,\beta)\geq 3 $ \emph{jointly fill up} $S$, i.e. are such
that every simple closed essential
curve $\zeta\in {\cal C}(S)$ intersects either
$\alpha$ or $\beta$ transversely,
by possibly increasing the lower bound for the diameter
of $J$ we may assume that
any two projective measured laminations $[\alpha],[\beta]$ which are
realized at the two distinct endpoints of $J$
jointly fill up $S$.

There is a 1-1-correspondence between measured
geodesic laminations and equivalence classes of
\emph{measured
foliations} on $S$ (see e.g. \cite{Ke} for a precise
statement and references). Via this identification,
any pair of distinct points
$[\lambda]\not=[\mu]\in {\cal P\cal M\cal L}$ which
jointly fill up the surface $S$ define
a unique Teichm\"uller geodesic
line. Thus
for every $\nu$-quasi-geodesic $\zeta:J\to {\cal C}(S)$
in ${\cal C}(S)$ with $\vert J\vert\geq 1/\epsilon_0$,
any pair of
projective measured laminations $
[\lambda],[\mu]$
realized at the two (possibly infinite)
endpoints of $\zeta$ defines
a unique Teichm\"uller geodesic $\eta([\lambda],[\mu])$.

Choose a number $R>2\chi$ and a smooth function
$\sigma:[0,\infty)\to [0,1]$ with $\sigma[0,\chi]\equiv 1$ and
$\sigma[R,\infty)\equiv 0$. For each $h\in {\cal T}_{g,m}$,
the number of simple closed curves
$\alpha$ on $S$ with $\ell_h(\alpha)\leq R$
is bounded from above by a universal constant not
depending on $h$, and the diameter of the subset
of ${\cal C}(S)$ containing these
curves is uniformly bounded as well.
Thus we obtain for every $h\in {\cal T}_{g,m}$ a finite Borel
measure $\mu_h$ on ${\cal C}(S)$ by defining
$\mu_h=\sum_\beta \sigma(\ell_h(\beta))\delta_\beta$
where $\delta_\beta$ denotes the Dirac mass at $\beta$.
The total mass of $\mu_h$ is bounded from
above and below by a universal positive constant, and
the diameter of the support of $\mu_h$ in ${\cal C}(S)$
is uniformly bounded as well. Moreover, the measures
$\mu_h$ depend continuously on $h\in {\cal T}_{g,m}$
in the weak$^*$-topology. This means that for every
bounded function $f:{\cal C}(S)\to \mathbb{R}$
the function $h\to \int f d\mu_h$ is continuous.

We define now a new ``distance'' function $\rho$ on
${\cal T}_{g,m}$ by
\[\rho(h,h^\prime)=\int_{{\cal C}(S)\times {\cal C}(S)}
d(\cdot,\cdot)d\mu_h\times d\mu_h^{\prime}/\mu_h({\cal C}(S))
\mu_{h^\prime}({\cal C}(S)).\]
Clearly the function $\rho$ is positive and
continuous on
${\cal T}_{g,m}\times {\cal T}_{g,m}$ and invariant
under the action of ${\cal M}_{g,m}$. Moreover, it
is immediate that there is a universal
constant $a>0$ such that $\rho(h,h^\prime)/a-a\leq
d(\Psi(h),\Psi(h^\prime))\leq a\rho(h,h^\prime)+a$.
As a consequence, for every $\nu>1$ there is a constant
$p=p(\nu)>1$ with the following property.
If $\gamma:J\to {\cal T}_{g,m}$ is such
that $\Psi\gamma$ is a $\nu$-quasi-geodesic, then
$\gamma$ is a $p$-quasi-geodesic with respect to
the ``distance'' function $\rho$. By this we mean that
\[\rho(\gamma(s),\gamma(t))/p-p\leq \vert s-t\vert
\leq p\rho(\gamma(s),\gamma(t))+p\]
for all $s,t\in J$. Moreover, for every $p>1$ there is a constant
$\nu=\nu(p)>1$ such that
if $\gamma:J\to {\cal T}_{g,m}$ is a Lipschitz curve which
is a $p$-quasi-geodesic with respect to $\rho$,
then $\Psi\circ
\gamma$ is a $\nu$-quasi-geodesic in ${\cal C}(S)$.

Let $h\in {\cal T}_{g,m}$ and let $\mu\in {\cal M\cal L}$
be a measured geodesic lamination.
The product of the transverse measure for $\mu$ together
with the length element of $h$ defines a measure on
the support of $\mu$
whose total mass is called the \emph{$h$-length} of $\mu$;
we denote it by $\ell_h(\mu)$.
Following Mosher \cite{Mo03}, for $p>1$ define
$\Gamma_{p}$ to be the set of all
triples
$(\gamma:J\to {\cal T}_{g,m},\lambda_+,\lambda_-)$ with the
following properties.
\begin{enumerate}
\item $0\in J$ and the diameter of
$J$ is at least $1/\epsilon_0$ where
$\epsilon_0=\epsilon_0(\nu(p))$ is as in Lemma 2.5.
\item $\gamma:J\to {\cal T}_{g,m}$ is a $p$-Lipschitz curve
which is a $p$-quasi-geodesic
with respect to the ``distance'' $\rho$.
\item $\lambda_+,\lambda_-\in {\cal M\cal L}$
are laminations of $\gamma(0)$-length 1, and the
projective measured lamination $[\lambda_+]$
is realized at
the right end, the projective measured lamination
$[\lambda_-]$
is realized at the left end
of $\gamma$.
\end{enumerate}

We equip $\Gamma_{p}$ with the product topology,
using the weak$^*$-topology on ${\cal M\cal L}$ for the
second and the third component of our triple and the compact-open
topology for the arc $\gamma$ in ${\cal T}_{g,m}$.
Note that this topology is metrizable.

We follow
Mosher (Proposition 3.17 of \cite{Mo03}) and show that
the action of ${\cal M}_{g,m}$ on $\Gamma_{p}$
is cocompact. Namely, recall from Lemma 2.5 that
there is a constant $\epsilon_0 >0$ such that for
every $(\gamma,\lambda_+,\lambda_-)\in \Gamma_p$ the image of
$\gamma$ is contained in ${\cal T}_{g,m}^{\epsilon_0}$.
Since ${\cal M}_{g,m}$ acts cocompactly on
${\cal T}_{g,m}^{\epsilon_0}$ it is therefore enough to show that
the subset of $\Gamma_p$ consisting of triples
with the additional property that $\gamma(0)$ is contained
in a fixed compact subset $A$ of ${\cal T}_{g,m}^{\epsilon_0}$ is
compact. Since our topology is metrizable, this follows if
every sequence of points $(\gamma,\lambda_+,\lambda_-)$
with $\gamma(0)\in A$ has a
convergent subsequence.

However, by the
Arzela-Ascoli theorem,
the set of $p$-Lipschitz maps into ${\cal T}_{g,m}^{\epsilon_0}$
issuing from a point in
$A$ is compact. Moreover, the function $\rho$ on
${\cal T}_{g,m}\times{\cal T}_{g,m}$ is continuous
and invariant under the action of ${\cal M}_{g,m}$ and
hence if $\gamma_i$ converges
locally uniformly to $\gamma$ and if $\gamma_i$
is a $p$-quasi-geodesic with respect to $\rho$
for all $i$ then the same is true for
$\gamma$. Since the function on ${\cal T}_{g,m}\times {\cal M\cal L}$
which assigns to a metric $h\in {\cal T}_{g,m}$ and a
measured lamination
$\mu\in {\cal M\cal L}$ the $h$-length of $\mu$ is continuous
and since for every fixed $h\in {\cal T}_{g,m}$ the set of
measured
laminations of $h$-length 1 is compact and
naturally homeomorphic to ${\cal P\cal M\cal L}$, the action of
${\cal M}_{g,m}$ on $\Gamma_p$ is indeed
cocompact provided that
the following holds: If $(\gamma_i:J_i\to {\cal T}_{g,m}^{\epsilon_0})$
is a sequence
of $p$-Lipschitz curves which converge locally uniformly
to $\gamma:J\to {\cal T}_{g,m}^{\epsilon_0}$,
if the projective measured lamination
$[\lambda_i]$ is realized at the right
endpoint of $J_i$ and if $[\lambda_i]\to [\lambda]$
in ${\cal P\cal M\cal L}$
$(i\to\infty)$ then
$[\lambda]$ is realized at the right endpoint of $J$.

To see that this is indeed the case, assume first
that $J\cap [0,\infty)=[0,b]$ for some $b\in (0,\infty)$.
Then for sufficiently large $i$ we have
$J_i\cap [0,\infty)=[0,b_i]$ with $b_i\in (0,\infty)$ and
$b_i\to b$. Thus $\gamma_i(b_i)\to \gamma(b)$
$(i\to\infty)$ and therefore for sufficiently large $i$
there is only a \emph{finite}
number of curves $\alpha\in {\cal C}(S)$ whose
length with respect to one of the metrics $\gamma_j(b_j),\gamma(b)$
$(j\geq i)$ is at most $\chi$. By passing
to a subsequence we may assume that there is a simple
closed curve $\alpha\in {\cal C}(S)$ with
$[\lambda_j]=[\alpha]$ for all large $j$.
The $\gamma_j(b_j)$-length of $\alpha$ is
at most $\chi$ for all sufficiently large $j$ and hence
the same is true for the $\gamma(b)$-length of $\alpha$
by continuity of the length function.
As a consequence, the limit $[\lambda]=[\alpha]$
of the sequence $([\lambda_i])$ is
realized at the endpoint $\gamma(b)$ of $\gamma$.

In the case that $[0,\infty)\subset J$ we argue as before.
Namely,
assume without loss of generality
that $b_i<\infty$ for all $i$ and that $b_i\to \infty$.
Recall that each of the curves $\Psi\gamma_i$ is
a uniform quasi-geodesic in ${\cal C}(S)$ and that
the map $\Psi$ is coarsely Lipschitz.
Let $\alpha_i\in {\cal C}(S)$ be the simple closed curve
such that $[\alpha_i]=[\lambda_i]$.
Then for each $i$, the curve $\alpha_i$ is contained
in a ball about $\Psi(\gamma_i(b_i))$ of radius
$R>0$ independent of $i$
and hence
as $i\to\infty$, the curves $\alpha_i$ converge
in ${\cal C}(S)\cup \partial {\cal C}(S)$
to the endpoint $\mu\in \partial {\cal C}(S)$
of $\Psi\circ\gamma$ in the Gromov
boundary of ${\cal C}(S)$. As a consequence, the curves
$\alpha_i$ converge to $\mu$ in the
\emph{coarse Haus\-dorff to\-po\-lo\-gy}
\cite{H04}. This means that every accumulation point
of $(\alpha_i)$ in the Hausdorff topology
contains $\mu$ as a sublamination.
Since $\mu$ is a minimal geodesic lamination
which fills up $S$, the complement of $\mu$ in
every lamination $\zeta$
containing $\mu$ as a sublamination
consists of a finite number of isolated leaves and
therefore every transverse measure supported
in $\zeta$ is in fact supported in $\mu$. Thus
after passing to a subsequence,
the projective measured laminations $[\lambda_i]$ converge
as $i\to \infty$ to a projective
measured lamination supported in $\mu$.
But $[\lambda_i]\to [\lambda]$ in ${\cal P\cal M\cal L}$ by
assumption and hence
the lamination $[\lambda]$ is realized
at the endpoint of $\gamma$. This shows our above claim
and implies that the action of ${\cal M}_{g,m}$ on $\Gamma_p$
is indeed cocompact.

Now we follow Section 3.10 of \cite{Mo03}.
Namely, each point $(\gamma,\lambda_-,\lambda_+)\in
\Gamma_p$
determines the geodesic $\eta([\lambda_-],[\lambda_+])$
in ${\cal T}_{g,m}$. This geodesic defines a family
$q_t$ of quadratic differentials whose horizontal foliation
corresponds to the lamination
$e^{-t}\lambda_+$ and whose vertical foliation
corresponds to $e^{t}\lambda_-$ (note that the
area of these
differentials is not necessarily normalized,
see \cite{Mo03}).

For $(\gamma,\lambda_-,\lambda_+)$ define $\sigma(\gamma,\lambda_-,
\lambda_+)$ to be the point on the geodesic
$\eta([\lambda_-],[\lambda_+])$ which
up to normalization corresponds
to the quadratic differential defined by
$\lambda_-,\lambda_+$.
The map taking $(\gamma,\lambda_-,\lambda_+)$
to $(\gamma(0),\sigma(\gamma,\lambda_-,\lambda_+))
\in {\cal T}_{g,m}\times
{\cal T}_{g,m}$ is continuous and equivariant with respect to
the natural action of ${\cal M}_{g,m}$
on $\Gamma_p$ and on ${\cal T}_{g,m}\times {\cal T}_{g,m}$.
Since the action
of ${\cal M}_{g,m}$ on $\Gamma_p$ is cocompact, the
same is true for the action of ${\cal M}_{g,m}$
on the image of our map (see \cite{Mo03}).
Thus the distance between $\gamma(0)$ and
$\sigma(\gamma,\lambda_-,\lambda_+)$
is bounded from above by a universal constant $R>0$.

Let again $(\gamma,\lambda_-,\lambda_+)\in \Gamma_p$.
For each $s\in J$ define
\[a_-(s)=\frac{1}{\ell_{\gamma(s)}(\lambda_-)},\quad a_+(s)=\frac{1}
{\ell_{\gamma(s)}(\lambda_+)}\]
where as before, $\ell_{\gamma(s)}(\lambda_{\pm})$ is the
$\gamma(s)$-length of $\lambda_{\pm}$.
These are continuous functions of $t\in J$.
Define for $s\in \mathbb{R}$ the shift $\gamma^\prime(t)=
\gamma(t+s)$; then the ordered triple
$(\gamma^\prime(0),a_-(s)\lambda_-,a_+(s)\lambda_+)$ lies
in the ${\cal M}_{g,m}$-cocompact set $\Gamma_{p}$ and hence
the distance between $\gamma(s)$ and a suitably chosen
point on the geodesic $\eta([\lambda_-],[\lambda_+])$ is
at most $R$. As a consequence, the arc $\gamma$ is
contained in the $R$-neighborhood of the geodesic
$\eta([\lambda_-],[\lambda_+])$.
Since the curve $\gamma$
is a $p$-quasi-geodesic,
this implies that the Hausdorff distance
between $\gamma(J)$ and a subarc
of $\eta([\lambda_-],[\lambda_+])$ connecting the
same endpoints is uniformly bounded and
shows the lemma.
\qed

\bigskip

Recall that for every finitely generated group $\Gamma$, every
finite symmetric set of generators induces a word norm on
$\Gamma$, and any two such word norms are equivalent.
As in the introduction, we define a convex
cocompact subgroup of the mapping class group
${\cal M}_{g,m}$ for $S$
as follows.

\bigskip

{\bf Definition 2.7:} Let $\Gamma$ be a finitely
generated subgroup of ${\cal M}_{g,m}$. The group
$\Gamma$ is called \emph{convex cocompact} if
some orbit map $\phi\in \Gamma\to \phi\alpha\in {\cal C}(S)$
for the action of $\Gamma$ on ${\cal C}(S)$
is a quasi-isometric embedding.

\bigskip

The following observation follows immediately
from the fact that the map $\Psi$ is coarsely
Lipschitz.

\bigskip

{\bf Lemma 2.8:} {\it For every
convex cocompact group $\Gamma<{\cal M}_{g,m}$ and every
$h\in {\cal T}_{g,m}$, the orbit map
$\phi\in \Gamma\to \phi h\in {\cal T}_{g,m}$
is a quasi-isometric embedding.}

{\it Proof:} Let $\Gamma<{\cal M}_{g,m}$ be a convex
cocompact group with a finite symmetric
set ${\cal G}$ of generators and let $h\in {\cal T}_{g,m}$.
Write $\ell=\max\{d(h,\phi h)\mid \phi\in {\cal G}\}$; since
$\Gamma$ acts on ${\cal T}_{g,m}$ as a group of isometries
we have $d(\phi h,\psi h)\leq \ell\Vert \phi^{-1}\psi\Vert$
for all $\phi,\psi\in \Gamma$. On the other hand,
the map $\Psi$ is coarsely equivariant with respect to
the action of ${\cal M}_{g,m}$ on ${\cal T}_{g,m}$
and ${\cal C}(S)$ and coarsely
Lipschitz; therefore there is a number
$\nu>0$ such that $d(\phi h,\psi h)\geq
d(\phi \Psi h,\psi \Psi h)/\nu-\nu$ for
all $\phi,\psi\in {\cal M}_{g,m}$. Since $\Gamma$
is convex cocompact by assumption, we moreover have
$d(\phi \Psi h,\psi \Psi h)\geq
\Vert \phi^{-1}\psi\Vert/\nu^\prime-\nu^\prime$
for some $\nu^\prime >0$ and all
$\phi,\psi\in\Gamma$. This shows the lemma.
\qed

\bigskip

The $\Gamma$-orbit $\Gamma h$
of a finite generated subgroup $\Gamma$ of ${\cal M}_{g,m}$
is called \emph{quasi-convex}
if for any two $\phi,\psi\in \Gamma$, the Teichm\"uller geodesic
connecting $\phi$ to $\psi$ is contained in a uniformly
bounded neighborhood of $\Gamma h$.
The following result shows the first part of
our theorem from the introduction in the more
general context of non-exceptional surfaces of finite type.
It was independently and at the same time
obtained by Kent and Leininger \cite{KL05}, with a different
proof.

\bigskip

{\bf Theorem 2.9:} {\it A finitely generated subgroup
$\Gamma$ of ${\cal M}_{g,m}$
is convex cocompact.
if and only if some $\Gamma$-orbit on ${\cal T}_{g,m}$ is quasi-convex.
}

{\it Proof:}
Let $\Gamma$ be a finitely generated convex cocompact
subgroup of ${\cal M}_{g,m}$.
Let $h\in {\cal T}_{g,m}^\epsilon$
for some $\epsilon >0$. By Lemma 2.8,
the orbit map $\phi\in \Gamma\to \phi h\in {\cal T}_{g,m}$ is
a quasi-isometric embedding. In particular, for any
two $\phi,\psi\in \Gamma$ the orbit of a geodesic
in $\Gamma$ connecting $\phi h$ to $\psi h$ is a uniform
quasi-geodesic in ${\cal T}_{g,m}$ which is contained in
${\cal T}_{g,m}^\epsilon$
and is mapped by $\Psi$ to a uniform quasi-geodesic in
${\cal C}(S)$.
By Theorem 2.1, this curve
is contained in a uniformly bounded neighborhood of the
Teichm\"uller geodesic connecting $\phi h$ to $\psi h$.
In other words, the orbit $\Gamma h\subset{\cal T}_{g,m}$
is quasi-convex.

Vice versa,
let $\Gamma<{\cal M}_{g,m}$ be a finitely generated group and
assume that there is some $h\in {\cal T}_{g,m}$ such that the
orbit $\Gamma h$ of $h$ for the action
of $\Gamma$ on ${\cal T}_{g,m}$ is quasi-convex. This means that
there is a number $D >0$ such that each Teichm\"uller geodesic
with both endpoints in $\Gamma h$ is contained in
the $D$-neighborhood of $\Gamma h$.
Using the map $\Psi:{\cal T}_{g,m}\to
{\cal C}(S)$ from the beginning of this Section,
write $\alpha=\Psi h\in {\cal C}(S)$
and assume to the contrary that the orbit map $\phi\to\phi\alpha$ for the
action of  $\Gamma$ on the complex of curves is \emph{not} a
quasi-isometry. Choose a finite symmetric set
${\cal G}$ of generators for $\Gamma$ which induces
the word norm $\Vert \,\Vert$.
Since $\Gamma$ acts on ${\cal C}(S)$ by
isometries, the orbit map
is coarsely Lipschitz with
respect to the word norm $\Vert\,\Vert$ on $\Gamma$
and the metric on ${\cal
C}(S)$. Thus our assumption implies that for every $L>0$ there is a
word $w=w_1\dots w_p\in \Gamma$ in the generators
$w_i\in {\cal
G}$  with $\Vert w\Vert =p$ and such that $d(w\alpha,\alpha)\leq
p/L$. Choose a geodesic $\zeta:[0,m]\to {\cal T}_{g,m}$ connecting
$h=\zeta(0)\in {\cal T}_{g,m}^\epsilon$ to $wh=\zeta(m)$. Since the
orbit $\Gamma h\subset {\cal T}_{g,m}$
is quasi-convex by assumption,
the geodesic is contained in the $D$- neighborhood of $\Gamma h$,
in particular it is contained in
${\cal T}_{g,m}^\epsilon$ for a universal number $\epsilon >0$. By
Theorem 2.1, there is a number $\nu>1$
not depending on $w$
such that the length $m$ of $\zeta$ is at most $\nu
p/L-1$.

On the other hand, the number of elements $\phi\in {\cal M}_{g,m}$
with $d(\phi h,h)\leq 2D+1$ is bounded from above by a constant
$\kappa>0$ only depending on the topological type of the surface
$S$ and on $\epsilon$ \cite {Bu}. Therefore
the word-norm of an element $\phi\in \Gamma$ with $d(\phi h,h)\leq
2D+1$ is bounded from above by a constant $\ell
>0$. For an integer
$k<m$ choose some $\phi(k)\in \Gamma$ with
$d(\phi(k)h,\zeta(k))\leq D$. Then we have
$d(\phi(k)^{-1}\phi(k+1)h,h)\leq 2D+1$ and hence the word norm of
$\phi(k)^{-1}\phi(k+1)$ is at most $\ell$. As a consequence, the
word norm $p$ of $w$ is at most $\ell (m+1)
\leq \ell \nu p/L$. For
$L>\ell\nu$, this is a contradiction which shows that $\Gamma$ is
indeed convex cocompact.
\qed

\bigskip

Improving earlier results of Farb and Mosher \cite{FM},
Kent and Leininger \cite{KL05}
obtained another characterization of convex
cocompact subgroups of ${\cal M}_{g,m}$. For the
formulation of their result,
let $\partial{\cal T}_{g,m}$ be the Thurston boundary
of Teichm\"uller space; recall that
$\partial {\cal T}_{g,m}$
is naturally homeomorphic to the space
${\cal P\cal M\cal L}$ of projective
measured laminations on $S$.
The \emph{limit set} $\Lambda$ for the action of
a group $\Gamma<{\cal M}_{g,m}$ on
Teichm\"uller space ${\cal T}_{g,m}$
is the closure in ${\cal P\cal M\cal L}$
of the set of fixed points
of the pseudo-Anosov elements of $\Gamma$; it is
a closed subset
of $\partial {\cal T}_{g,m}\sim
{\cal P\cal M\cal L}$ which is invariant under
the natural action of $\Gamma$ (see
\cite{MP89} for the construction of
limit sets for a subgroup of ${\cal M}_{g,m}$ and their basic
properties).
Its \emph{weak convex hull} in ${\cal T}_{g,m}$ is defined to
be the closure of the union of all geodesics in ${\cal T}_{g,m}$
with both endpoints in $\Lambda$; it is a closed
$\Gamma$-invariant subset of ${\cal T}_{g,m}$.

\bigskip

{\bf Theorem 2.10 \cite{KL05}:} {\it
A finitely generated subgroup $\Gamma$
of ${\cal M}_{g,m}$
is convex cocompact if and only if
the group $\Gamma$ acts cocompactly on the weak
convex hull in ${\cal T}_{g,m}$
of its limit set $\Lambda\subset {\cal P\cal M\cal L}$.}

\bigskip

We refer to \cite{KL05} for additional characterizations
of convex cocompact subgroups of ${\cal M}_{g,m}$ which
correspond to the various characterizations
of convex cocompact Kleinian groups.

\section{Hyperbolic $\mathbb{R}$-bundles over
proper hyperbolic spaces}

In \cite{FM}, Farb and Mosher introduce
\emph{metric fibrations} as a generalization of
Riemannian submersions between complete Riemannian
manifolds. The following definition is adapted to our
needs from their paper.

\bigskip

{\bf Definition 3.1:} Let $(X,d)$ be a proper geodesic metric
space. A \emph{metric fibration} over $X$ with fibre a topological
space $F$ is a geodesic metric space $(Y=X\times F,d)$ with the
following properties.
\begin{enumerate}
\item For all $x,x^\prime\in X$ and
every $y\in F$ we have $d((x,y),(x^\prime,y))=
d(x,x^\prime)=d(\{x\}\times F,\{x^\prime\}\times F)$.
\item For each $x\in X$, the metric
on $Y$ induces a complete geodesic metric on $\{x\}\times F$ which
defines the given topology on $\{x\}\times F\sim F$.
\end{enumerate}
The metric fibration is called \emph{bounded} if
there is a number $n>0$ such that for all $x,x^\prime\in X$,
the natural map $\{x\}\times F\to \{x^\prime\}\times F$
is bilipschitz with Lipschitz constant bounded from above
by $e^{nd(x,x^\prime)}$.

\bigskip

Recall that a geodesic metric space $(X,d)$ is called
\emph{$\delta$-hyperbolic} for some $\delta \geq 0$ if the
\emph{$\delta$-thin triangle condition} holds for $X$: For every
geodesic triangle in $X$ with sides $a,b,c$, the side $a$ is
contained in the $\delta$-neighborhood of $b\cup c$. In this
section we consider a metric fibration $Y=X\times T\to X$ over a
$\delta$-hyperbolic geodesic metric space $(X,d)$ with fibre a
simplicial tree $T$ of bounded valence. Our goal is to give a
necessary and sufficient condition for the space $Y$ to be
hyperbolic.

We begin with analyzing the case when $T$ is
a closed subset of the real line $\mathbb{R}$.
We use an idea of Bestvina and Feighn
who introduced in \cite{BF} the following ``rectangle flare''
condition. Let $(X,d)$ be a geodesic metric space and let $r:X\to
(0,\infty)$ be a positive function. Given $\kappa >1, n\in
\mathbb{Z}_+$, the \emph{$\kappa,n$-flaring property} for $r$ with
threshold $A\geq 0$ says that if $J\subset \mathbb{R}$ is a closed
connected subset, if $t-n,t,t+n\in J$ and if $\gamma:J\to X$ is a
geodesic so that $r(\gamma(t))\geq A$ then
$\max\{r(\gamma(t-n)),r(\gamma(t+n))\}\geq \kappa r(t)$.
We say
that $r$ satisfies the \emph{bounded} $\kappa,n$-flaring property
with threshold $A$ if in addition the growth of $r$ is uniformly
exponentially bounded with exponent $n$, i.e. if $r(y)\leq
e^{d(x,y) n}r(x)+nd(x,y)$ for all $x,y\in X$.
Note that this implies in particular that the function $r$ is
continuous.

For a constant $c>0$ define a subset $B$ of $X$
to be \emph{$c$-quasi-convex} if every geodesic in
$X$ connecting two points in $B$ is contained in the
$c$-neighborhood of $B$.
We have.

\bigskip

{\bf Lemma 3.2:} {\it Let $X$ be a proper
geodesic metric space and let $r:X\to (0,\infty)$ be a function
which satisfies the bounded
$\kappa,n$-flaring property with threshold
$A>0$. Let $\mu=\inf_{x\in X}r(x)$. There is a constant
$D=D(\kappa,n)>0$ only depending on
$\kappa,n$ with the following properties.
\begin{enumerate}
\item If $\mu\geq A$ then the
function $r$ assumes a minimum on
$X$. The diameter of the set
$\{x\in X\mid r(x)=\mu\}$ is bounded from above by $D$.
\item
If $\mu< A$ then the set
$\{x\mid r(x)\leq A\}$ is $D$-quasi-convex.
\end{enumerate}
}

{\it Proof:} Let $(X,d)$ be a proper geodesic metric space and let
$r:X\to (0,\infty)$ be a function which satisfies the bounded
$\kappa,n$-flaring property with threshold $A$. Assume that
$\mu=\inf_{x\in X}r(x)\geq A$. We have to show that $r$ assumes a
minimum on $X$. Namely, let $x,y$ be two points in $X$ whose
distance $\chi$ is at least $2n$ and such that $r(x)<\kappa \mu,
r(y)<\kappa \mu$. Let $\gamma:[0,\chi]\to X$ be a geodesic
connecting $\gamma(0)=x$ to $\gamma(\chi)=y$ and let $\ell\geq 2$
be such that $\chi\in [\ell n,(\ell +1)n)$. By the flaring
property for $r$ and the fact that $r(\gamma(n))\geq A$ we have
$r(\gamma(2n))\geq \kappa r(\gamma(n))\geq\kappa \mu$ and
inductively we conclude that $r(\gamma(\ell n))\geq \kappa^{\ell
-1}\mu$. On the other hand, the growth of $r$ is uniformly
exponentially bounded and therefore $\kappa \mu >r(y)\geq
e^{-n^2}r(\gamma(\ell n))-n^2\geq e^{-n^2}\kappa^{\ell
-1}\mu-n^2$. This implies that the distance between $x$ and $y$ is
bounded from above by a constant $D>0$ only depending on
$\kappa,n$. Since $X$ is proper and $r$ is continuous, we conclude
that the function $r$ assumes a minimum, and the diameter of the
set of points at which such a minimum is achieved is at most $D$.

Now assume that $\mu< A$ and let $E=\{z\mid r(z)\leq A\}$. We have
to show that $E$ is $D^\prime$-quasi-convex for a constant
$D^\prime >0$ only depending on $\kappa,n$. For this let $x,y\in
E$ and let $\gamma:[0,\chi]\to X$ be a geodesic arc connecting $x$
to $y$. Let $\ell\geq 0$ be such that the length $\chi$ of
$\gamma$ is contained in the interval $[\ell n,(\ell+1)n)$. If
$\ell\leq 1$ then there is nothing to show, so assume otherwise.
If $\gamma(n)\not\in E$ then we have $r(\gamma(n))> A$ and it
follows as above from the flaring property that $r(y)\geq
e^{-n^2}\kappa^{\ell-1}A-n^2$. Hence the distance $\chi$ between
$x$ and $y$ is bounded from above by a universal constant.
Otherwise we have $\gamma(n)\in E$ and we can apply the same
consideration to the points $\gamma(n),y$. Inductively we conclude
that the set $E$ is $D^\prime$-quasi-convex for a constant
$D^\prime >0$ only depending on $\kappa,n$. \qed

\bigskip

In the sequel, we will use the following criterion for
hyperbolicity of a geodesic metric space (Proposition 3.5 in
[H05]).

\bigskip

{\bf Lemma 3.3:} {\it Let $(Y,d)$ be a geodesic metric space.
Assume that
there is a number $D>0$ and for every pair of points
$x,y\in Y$ there is a path $c(x,y):[0,1]\to Y$ connecting
$c(x,y)(0)=x$ to $c(x,y)(1)=y$ with the following properties.
\begin{enumerate}
\item If $d(x,y)\leq 1$ then the diameter of the set
$c(x,y)[0,1]$ is at most $D$.
\item For $x,y\in X$ and $0\leq s\leq t\leq 1$, the
Hausdorff distance between $c(x,y)[s,t]$ and
$c(c(x,y)(s),c(x,y)(t))[0,1]$ is at most $D$.
\item For any triple $(x,y,z)$ of
points in $X$,
the arc $c(x,y)[0,1]$ is
contained in the $D$-neighborhood of $c(x,z)[0,1]\cup
c(z,y)[0,1]$.
\end{enumerate}
Then the space $(Y,d)$ is $\delta$-hyperbolic for a constant
$\delta >0$ only depending on $D$.
Moreover, for all $x,y\in Y$ the Hausdorff distance
between $c(x,y)$ and a geodesic connecting $x$ to $y$
is at most $\delta$.
}

\bigskip

Now consider a metric fibration whose fibre $J$ either is the
closed interval $[0,1]$ or the half-line $[0,\infty)$. By
assumption, for every compact interval $[s,t]\subset J$ and every $x\in X$
the arc $\{x\}\times [s,t]$ is rectifiable. As a consequence, we
can define a function on $X$ by associating to $x\in X$ the length
of the arc $\{x\}\times [s,t]$; we call such a function a
\emph{vertical distance function}. The next lemma is the main
technical result of this section.

\bigskip

{\bf Lemma 3.4:} {\it Let $(X,d)$ be a proper $\delta$-hyperbolic
geodesic metric space and let $Y=X\times J\to X$ be a
bounded metric
fibration. Assume that the vertical distance functions satisfy the
$\kappa,n$-flaring property with flaring threshold $A$ for
some $\kappa >1,n>0,A>0$. Assume moreover that the infimum of
every vertical distance function is not bigger than $A$. Then $Y$
is $\delta_0$-hyperbolic for a number $\delta_0 >0$ only depending
on $\delta,\kappa,n,A$.}

{\it Proof:} Let $(X,d)$ be a proper $\delta$-hyperbolic geodesic
metric space and let $Y=X\times J\to X$ be a bounded metric
fibration with fibre $J=[0,1]$ or $J=[0,\infty)$ such that the
vertical distance functions
satisfy the $\kappa,n$-flaring property with
threshold $A>0$ for some $\kappa >1,n>0,A>0$.
By the definition of a bounded metric fibration,
the vertical distance functions satisfy in fact the \emph{bounded}
$\kappa,n$-flaring property with threshold $A$. Denote the
distance on $Y$ again by $d$. For $t\in J$ let $\ell_t:X\to
[0,\infty)$ be the function which associates to a point $x\in X$
the length of the vertical path $\{x\}\times [0,t]$. By
assumption, the function $\ell_t$ is continuous. Write
$\mu(t)=\inf_{x\in X} \ell_t(x)$; the function $t\to \mu(t)$ is
continuous and increasing. We assume that $\mu$ is bounded from
above by $A$. Our goal is to construct for any two points $x,y\in
Y$ a curve $c(x,y)$ connecting $x$ to $y$ so that the resulting
curve system satisfies the properties 1-3 in Lemma 3.3. For this
we proceed in four steps.

{\sl Step 1:}

In a first step, we construct for every
$y=(x,t)\in Y$ and every $s\leq t$ a curve
$\eta_s(x,t):[0,1]\to Y$ connecting
$(x,t)$ to $X\times \{s\}$. For this
let $\overline{X}$
be the union of $X$ with its Gromov boundary $\partial X$. Since
$X$ is proper, the space $\overline{X}$ is compact.
For $s\geq 0$ define a set
$C_s\subset \overline{X}$ as follows. If $J=[0,1]$ then define
$C_s=\{x\in X\mid (\ell_1-\ell_s)(x)\leq A\}$.
By Lemma 3.2, the set $C_s$ is
$D_1$-quasi-convex for a universal constant $D_1>0$. If
$J=[0,\infty)$ then for $t\geq s$ write $Q_{t,s}=
\{\ell_t-\ell_s\leq A\}$.
By Lemma 3.2 the sets $Q_{t,s}$ are $D_1$-quasi-convex for our
constant $D_1>0$.
If we denote by
$\overline{Q}_{t,s}$ the closure of $Q_{t,s}$ in $\overline{X}$
then the sets $\overline{Q}_{t,s}$ are compact and non-empty
and we have
$\overline{Q}_{t,s}\supset \overline{Q}_{u,s}$ for $t\leq u$.
Thus $C_s=\cap_{t\geq s}\overline{Q}_{t,s}\not=\emptyset$,
moreover $C_s$ is $D_1$-quasi-convex. This means
that either $C_s$ consists of a single point
$\zeta\in \partial X$ or $C_s\cap X$ is
non-empty and $D_1$-quasi-convex, with closure
$C_s$ in $\overline{X}$.

For $(x,t)\in Y$ and $s\leq t$
define now a curve $\eta_s(x,t):[0,1]\to Y$ connecting
$(x,t)=\eta_s(x,t)(0)$ to $\eta_s(x,t)(1)\in X\times \{s\}$
as follows. First, if $t=s$
then let $\eta_s(x,t)(\tau)=(x,t)$ for all
$\tau\in [0,1]$. If $y=(x,t)$ for some $t>s$ then choose a minimal
geodesic $\gamma_{x,s}:[0,\sigma]\to X$ connecting
$\gamma_{x,s}(0)=x$ to the set
$C_s$; if $C_s$ is a single point $\zeta\in
\partial X$ then $\sigma=\infty$ and we require that
$\gamma_{x,s}$ converges to $\zeta$. We assume that the choice of
$\gamma_{x,s}$ only depends on $x,s$ but not on $t$. Since
$\ell_t<\ell_{t^\prime}$ for $t< t^\prime$, there is a smallest
number $\nu_{x,t,s}\geq 0$ so that
$(\ell_t-\ell_s)(\gamma_{x,s}(\nu_{x,t,s}))\leq A$. Let
$\eta_s(x,t)$ be a reparametrization on the interval $[0,1]$ of
the horizontal arc $\gamma_{x,s}[0,\nu_{x,t,s}]\times \{t\}$ with
a vertical arc of length at most $A$ connecting
$(\gamma_{x,s}(\nu_{x,t,s}),t)$ to
$(\gamma_{x,s}(\nu_{x,t,s}),s)\in X\times \{s\}$.

{\sl Step 2:}

In a second step, we show that for every
$R>0$, $s\in J$ and all $y,z\in X\times ([s,\infty)\cap J)$
with $d(y,z)\leq R$ the Hausdorff distance between
$\eta_s(y)$ and $\eta_s(z)$ is bounded from
above by a number $\tau(R)>0$ only depending on
$R$ but not on $s,y,z$.

We first consider the case that the points $y,z$ are contained in
$X\times \{t\}$ for some fixed $t\geq s$. Thus let $R>0$, let
$t\geq 0$, let $x,u\in X$ with $d(x,u)\leq R$, let
$y=(x,t),z=(u,t)\in Y$ and let $s\in [0,t]$. By hyperbolicity of
$X$, the Hausdorff distance between the two geodesics
$\gamma_{x,s},\gamma_{u,s}$ of minimal length connecting the
points $x,u$ of distance at most $R$ to the $D_1$-quasi-convex
subset $C_s$ of $\overline{X}$ is bounded from above by a
universal constant $\tau_1(R)>0$ only depending on $R$ but not on $x,u$.

There are smallest
numbers $\nu_{x,t,s}\geq 0,\nu_{u,t,s}\geq 0$ such that
$(\ell_t-\ell_s)(\gamma_{v,s}(\nu_{v,t,s}))\leq A$ $(v=x,u)$.
By the definition of the curves $\eta_s(y)$ it is now enough
to show that the distance between
$\gamma_{x,s}(\nu_{x,t,s})$ and $\gamma_{u,s}(\nu_{u,t,s})$
is bounded from above by a constant which only
depends on $R$.
Note that for $s=t$ we have
$\nu_{x,t,s}=0=\nu_{u,t,s}$ and hence there is nothing to show, so
assume that $s<t$.

Since by assumption the growth of the vertical distance
functions is uniformly exponentially bounded, there is a universal
number $\beta>0$ only depending on $\tau_1(R)$ such that for any
two points $v,w\in X$ with $d(v,w)\leq \tau_1(R)$ we have
$(\ell_t-\ell_s)(v)\leq \beta(\ell_t-\ell_s)(w)$.
By the definition of $\nu_{x,t,s}$ and the
flaring property, this means that there is a number $\xi >0$ only
depending on $R$ such that $(\ell_t-\ell_s)(w)>A$ whenever $w\in
X$ is such that $d(w,\gamma_{x,s}(\sigma))\leq \tau_1(R)$
for some $\sigma\in
[0,\nu_{x,t,s}-\xi]$ (compare the argument in the proof of Lemma 3.2). As
a consequence, if $a\geq 0$ is such that
$d(\gamma_{x,s}(a),\gamma_{u,s}(\nu_{u,t,s}))
\leq \tau_1(R)$ then $a\geq
\nu_{x,t,s}-\xi$. Since $\gamma_{x,s},\gamma_{u,s}$ are
geodesics of Hausdorff distance at most $\tau_1(R)$ we conclude
that $d(\gamma_{u,s}(0),\gamma_{u,s}(\nu_{u,t,s})) \geq
d(\gamma_{x,s}(0),\gamma_{x,s}(a))-R-\tau_1(R)=a-R-\tau_1(R)
\geq \nu_{x,t,s}-\xi-R-\tau_1(R)$.
Exchanging the role of $x$ and $u$ then shows that
$\vert \nu_{x,t,s}-\nu_{u,t,s}\vert
\leq \xi+R+\tau_1(R)$. But $\gamma_{x,s},\gamma_{u,s}$
are geodesics of Hausdorff distance at most $\tau_1(R)$ and
therefore the distance
between $\gamma_{x,s}(\nu_{x,t,s})$ and $\gamma_{u,s}(\nu_{u,t,s})$
is indeed bounded from
above by a universal constant $\tau_2(R)>0$ only depending on $R$.
As a consequence, the Hausdorff distance between $\eta_s(y)$ and
$\eta_s(z)$ is bounded from above by a number $\tau_3(R)>0$ only
depending on $R$.

Now consider nearby
points $y,z$ contained in the same fibre of our metric fibration.
Thus let $x\in X$ and let $t\geq 0,b>0$ be such that the length of
the vertical arc $\{x\}\times [t,t+b]$ is at most $A$. Write
$y=(x,t),z=(x,t+b)$ and let $s\leq t$;
we claim that the Hausdorff distance between
$\eta_s(y)$ and $\eta_s(z)$ is bounded
from above by a universal constant.

Namely, let again
$\gamma_{x,s}$ be the minimal geodesic
connecting $x$ to $C_s$ as in the definition of the arc
$\eta_s(x,t)$.
There is a minimal number $\sigma_b=\nu_{x,t+b,s}\geq 0$ such that
$(\ell_{t+b}-\ell_s)(\gamma_{x,s}(\sigma_b))\leq A$
and hence
$\ell_{t+b}(\gamma_{x,s}(\sigma))-
\ell_t(\gamma_{x,s}(\sigma))\leq A$ for $\sigma=0,\sigma_b$.
By the bounded $\kappa,n$- flaring property
for vertical distances, this
implies that there is a universal number $A^\prime \geq A$ such
that $(\ell_{t+b}-\ell_t)(\gamma_{x,s}(\sigma))\leq A^\prime$
for all $\sigma\in [0,\sigma_{b}]$. Let
$\sigma_0=\nu_{x,t,s}\leq \sigma_b$
be the minimal number such
that $(\ell_t-\ell_s)(\gamma_{x,s}(\sigma_0))\leq A$.
Then
$A\leq (\ell_{t+b}-\ell_s)(\gamma_{x,s}(\sigma))\leq 2A^\prime$ for every
$\sigma\in [\sigma_0,\sigma_b]$
and consequently an application of the flaring property
as in the proof of Lemma 3.2
for the function $\ell_{t+b}-\ell_s$ yields that
$\sigma_{b}\leq \sigma_0+\xi$
for a universal number $\xi >0$. It follows that
the Hausdorff distance between
$\eta_s(y)$ and $\eta_s(z)$
is bounded from above by a universal constant.

By our assumption that the growth of the vertical
distance functions is uniformly exponentially bounded, for
every $R>0$ there is a number $\nu(R)>0$ such that
for every $y=(x,t)\in Y$ the $R$-ball about
$y$ is contained in the set $\{z=(u,s)\in Y\mid
d(x,u)\leq \nu(R),\vert (\ell_t-\ell_s)(y)\vert \leq \nu(R)\}$.
Together we conclude that for every $R>0$ we can
find a number $\tau(R)>0$ with the following property.
Let $y=(x,t),y^\prime=(x^\prime,t^\prime)\in Y$ with
$d(y,y^\prime)\leq R$; then for every $s\leq \min\{t,t^\prime\}$
the Hausdorff distance between $\eta_s(y)$ and
$\eta_s(z)$ is bounded from above by $\tau(R)$.

{\sl Step 3:}

Define a system $c(y,z)$ of arcs connecting an arbitrary pair of
points $y,z\in Y$ as follows. If $y=(x,t),z=(u,s)\in Y$ with
$0\leq s\leq  t$ then define $c(y,z)$ to be a reparametrization on
$[0,1]$ of the composition of the arc $\eta_s(y)$ with a geodesic
in $X\times \{s\}\sim X$ connecting $\eta_s(y)(1)$ to $z$. Define
also $c(z,y)$ to be the inverse of $c(y,z)$.

In our third step we
show that for every $R>0$ and all
$y,y^\prime\in Y$ with
$d(y,y^\prime)\leq R$, all $z\in Y$  the Hausdorff distance
between $c(y,z),c(y^\prime,z)$ is bounded from
above by a constant $\chi(R)>0$ only depending on $R$.
For this let
$R>0$ and let
$y,y^\prime,z\in Y$ with $d(y,y^\prime)\leq R$.
We distinguish 3 cases.

{\sl Case 1:} $z=(u,s),y=(x,t),y^\prime=(x^\prime,t^\prime)$ with
$0\leq s\leq t\leq t^\prime$.

By the definition of the curves $c(y,z)$ and Step 2, the curves
$c(y,z),c(y^\prime,z)$ are composed of the arcs $\eta_s(y),
\eta_s(y^\prime)$ of Hausdorff distance at most $\tau(R)$ and
geodesic arcs in $X\times \{s\}\sim X$ connecting the points
$\eta_s(y)(1),\eta_s(y^\prime)(1)$ of distance at most $\tau(R)$
to $z$. By $\delta$-hyperbolicity of $X\times \{s\}\sim X$, the
Hausdorff distance between $c(y,z)$ and $c(y^\prime,z)$ is bounded
from above by a constant $\chi_1(R)>0$ only depending on $R$.

{\sl Case 2:} $z=(u,s),y=(x,t),y^\prime=(x^\prime,t^\prime)$ with
$0\leq t\leq s\leq t^\prime$.

Since the distance between $y$ and
$y^\prime$ is at most $R$ and $Y$ is a geodesic metric space,
there is a point $y^{\prime\prime}=(x^{\prime\prime},s)\in X\times
\{s\}$ whose distance to both $y,y^\prime$ is at most $R$. By Case
1 above, the Hausdorff distance between $c(y^\prime,z)$ and
$c(y^{\prime\prime},z)$ is at most $\chi_1(R)$. Thus we may assume
without loss of generality that $t^\prime=s$; then $c(y^\prime,z)=
c(z,y^\prime)$ is the lift to $X\times \{s\}$ of a geodesic in $X$
connecting $u$ to $x^\prime$. Since $d(y,y^{\prime})\leq R$ and
$y^\prime=(x^\prime,s),y=(x,t)$, we have $d(x,x^\prime)\leq R$ and
hence by hyperbolicity of $X$, the Hausdorff distance between a
geodesic connecting $u$ to $x$ and a geodesic connecting $u$ to
$x^\prime$ is bounded by a uniform constant $\chi_2(R)>0$. Thus
the Hausdorff distance between $c(y^\prime,z)$ and $c((x,s),z)$ is
at most $\chi_2(R)$ and we may assume without loss of generality
that $x=x^\prime$.

By the flaring property for vertical distances, the point $x$ is
contained in a uniformly bounded neighborhood of the set
$E=\{\ell_s-\ell_t\leq A\}$. Since $E$ is $D_1$-quasi-convex, by
hyperbolicity a geodesic connecting $u$ to $x$ is contained in a
uniform neighborhood of the composition of a minimal geodesic
$\zeta$ connecting $u$ to $E$ and a geodesic arc connecting the
endpoint of $\zeta$ to $x$. By construction, the curve
$c(y,z)=c(z,y)$ is composed of the lift to $X\times \{s\}$ of a
minimal geodesic $\zeta:[0,\tau]\to X$ connecting $u$ to $E$, a
vertical arc of length at most $A$ connecting $(\zeta(\tau),s)$ to
$(\zeta(\tau),t)$ and the lift to $X\times \{t\}$ 
of a geodesic $\xi$ in
$X$ connecting $\zeta(\tau)$ to $x$ which is contained in a
uniformly bounded neighborhood of $E$. It follows that the
Hausdorff distance between $c(y,z)$ and the lift of the
composition of $\zeta$ and $\xi$ to $X\times \{t\}$ is uniformly
bounded. Therefore the Hausdorff distance between
$c(y,z),c(y^\prime,z)$ is bounded from above by a constant
$\chi_3(R)>0$ only depending on $R$.

{\sl Case 3:} $z=(u,s),y=(x,t),y^\prime=(x^\prime,t^\prime)$ with
$0\leq t\leq t^\prime\leq s$.

We claim that $\eta_t(z)$ contains a subarc whose Hausdorff
distance to $\eta_{t^\prime}(z)$ is uniformly bounded. Namely, the
sets $C_t,C_{t^\prime}\subset X$ are $D_1$-quasi-convex and
$C_t\subset C_{t^\prime}$. By hyperbolicity of $X$, if $u\in
X-C_{t^\prime}$ then a minimal geodesic $\xi$ in $X$ connecting
$u$ to $C_t$ is contained in a uniformly bounded neighborhood of
the composition of a minimal geodesic $\zeta_1:[0,a]\to X$
connecting $u$ to $C_{t^\prime}$ and a minimal geodesic
$\zeta_2$ connecting $\zeta_1(a)=\zeta_2(0)$ to $C_t$. From this
and the definition of our arcs $\eta_t$, the claim is immediate.

Denote by $z^{\prime\prime}=\eta_t(z)(\sigma)$ $(\sigma\in [0,1])$
the endpoint of this subarc and write
$z^{\prime}=\eta_{t^\prime}(z)(1)\in X\times \{t^\prime\}$. The
curve $c(y^\prime,z)$ is composed of the arcs $\eta_{t^\prime}(z)$
and $c(y^\prime,z^\prime)$,
and the curve $c(y,z)$ is composed of the arcs 
$\eta_t(z)[0,a]$ and $c(y,z^{\prime\prime})$.
Thus up to a constant only depending on $R$, the Hausdorff
distance between $c(y,z)$ and $c(y^\prime,z)$ is bounded from
above by the Hausdorff distance between $c(y,z^{\prime\prime})$ and
$c(y^\prime,z^\prime)$. Now the distance between $z^\prime$ and
$z^{\prime\prime}$ is uniformly bounded and hence by Case 1 above,
the Hausdorff distance between $c(y,z^\prime)$ and
$c(y,z^{\prime\prime})$ is bounded by a constant only depending on
$R$. In other words, for our estimate we may replace $z$ by
$z^\prime$, i.e. we may assume without loss of generality that
$s=t^\prime$. However, this case is contained in Case 2 above.

Together we established an upper bound $\chi(R)>0$ for the
Hausdorff distance between $c(y,z)$ and $c(y^\prime,z)$ whenever
$d(y,y^\prime)\leq R$.

{\sl Step 4:}

In a final step, we show that
our system of curves satisfies the properties 1)-3) in Lemma 3.3.

Namely, for $y=z$ the curve $c(y,z)$ is constant, and hence if
$d(y,z)\leq 1$ then the diameter of $c(y,z)$ is at most $\chi(R)$
where $\chi(R)>0$ is as in Step 3. This means that property 1 is
valid with $D=\chi(1)>0$.

Similarly, let $y,z\in Y$ and let $0\leq s\leq t\leq 1$. Then
either the restriction of the curve $c(y,z)$ to $[s,t]$ is
obtained from our above procedure, i.e. from the same construction
used for the curve $c(c(y,z)(s),c(y,z)(t))$, or one of the points
$c(y,z)(s),c(y,z)(t)$ is contained in the vertical subarc of
$c(y,z)$. In the first case it is immediate from Step 3 above that
the Hausdorff distance between $c(y,z)[s,t]$ and
$c(c(y,z)(s),c(y,z)(t))$ is uniformly bounded. In the second case,
if say the point $c(y,z)(s)$ is contained in the vertical subarc
of $c(y,z)$ then there is some $s^\prime\geq s$ such that
$c(y,z)[s^\prime,t]$ is obtained from the above procedure and that
the Hausdorff distance between $c(y,z)[s,t]$ and
$c(y,z)[s^\prime,t]$ is bounded from above by a universal
constant. By Step 3, the Hausdorff distance between
$c(c(y,z)(s),c(y,z)(t))$ and $c(y,z)[s^\prime,t]$ is bounded from
above by a universal constant as well. As a consequence, there is
a number $\nu>0$ such that property 2 is valid with $D=\nu$.

We are left with showing that the $\delta_0$-thin triangle
condition for a universal number $\delta_0>0$ also holds. For this
let $y_1,y_2,y_3$ be any 3 points in $Y$. Assume that
$y_i=(x_i,s_i)$ with $0\leq s_1\leq s_2\leq s_3$. By construction
and Step 3 above, the curves $c(y_1,y_3),c(y_2,y_3)$ both contain
a subarc whose Hausdorff distance to $\eta_{s_2}(y_3)$ is
uniformly bounded. By our above estimate of Hausdorff distances,
this means that for the purpose of establishing the thin triangle
condition we may replace $y_3$ by $\eta_{s_2}(y_3)(1)$, i.e. we
may assume that in fact $s_2=s_3=s$. Then the arc $c(y_2,y_3)$ is
the lift to $X\times \{s\}$ of a geodesic $\gamma$ in $X$
connecting $x_2$ to $x_3$.

Let $E=\{u\in X\mid (\ell_{s}-\ell_{s_1})(u)\leq A\}$. Recall that
$E$ is $D_1$-quasi-convex. Let
$\zeta_i:[0,\sigma_i]\to X$ $(i=2,3)$ be a
minimal geodesic connecting $x_2,x_3$ to $E$. Then $\gamma$ is
contained in a uniformly bounded neighborhood of
the union of
$\zeta_2[0,\sigma_2]\cup \zeta_3[0,\sigma_3]$ with a geodesic arc
connecting $\zeta_2(\sigma_2)$ to $\zeta_3(\sigma_3)$.
Moreover, by our above considerations
the curves $c(y_2,y_1)$ and $c(y_3,y_1)$ contain
each a subarc $\nu_2,\nu_3\subset X\times \{s\}$ whose
Hausdorff distance to the arcs $\zeta_2\times \{s\},
\zeta_3\times \{s\}$ is uniformly bounded. As a consequence,
for the purpose of the thin triangle condition we may
as well assume that the points $y_2,y_3$ are contained
in a uniformly bounded neighborhood of $E$. However,
in this case the thin triangle condition is immediate
from the definition of the curves $c(x,y)$ and hyperbolicity of
$X$. As a consequence, our system of curves
$c(x,y)$ satisfies the properties 1)-3) in Lemma 3.3 for a number $D>0$
only depending on $\delta,\kappa,n,A$ and hence
the space $Y$ is $\delta^\prime$-hyperbolic for a
constant $\delta^\prime>0$ only depending on $\delta,\kappa,n,A$.
\qed

\bigskip

Let $X$ be a proper $\delta$-hyperbolic geodesic metric space.
Recall that a closed subset $E$ of $X$ is \emph{strictly convex}
if every geodesic connecting two points in $E$ is contained in
$E$. The following lemma shows that under suitable assumptions,
hyperbolicity is preserved under glueing along stricly convex
subsets. For its formulation, for a number $R>0$ we call two
closed strictly convex subsets $D,E$ of a $X$ \emph{$R$-separated}
if $D,E$ are disjoint and if moreover the following holds. Let
$\gamma:[0,a]\to X$ be a minimal geodesic connecting $D$ to $E$;
then $\gamma[0,a]$ is contained in the $R$-neighborhood of
\emph{every} geodesic connecting $D$ to $E$. For example, two
non-intersecting geodesics in the hyperbolic plane are
$R$-separated for a constant $R>0$ which tends to infinity as the
distance between the geodesics tends to zero. The two boundary
geodesics of a flat strip in $\mathbb{R}^2$ are not $R$-separated
for any $R>0$. Note also that by the explicit construction of the
curves $c(x,y)$ in the proof of Lemma 3.4 the following holds. If
$Y=X\times [0,1]\to X$ is a bounded metric fibration over a
hyperbolic geodesic metric space such that the vertical distance
functions satisfy the $\kappa,n$-flaring property with threshold
$A$ for some $\kappa >1,n>0,A>0$ and if the infimum of the
vertical length of the fibres equals $A$ then $Y$ is hyperbolic
and the subsets $X\times \{0\},X\times \{1\}$ are strictly convex
and $R$-separated for a number $R>0$ only depending on
$\kappa,n,A$.

\bigskip

{\bf Lemma 3.5:} {\it Let $\delta >0,R>0$,
let $I\subset\mathbb{Z}$
be any subset and let $X$ be
a geodesic metric space with the following
properties.
\begin{enumerate}
\item[a)] $X=\cup_{i\in I}X_i$ where for each $i\in I$,
$X_i$ is a proper $\delta$-hyperbolic geodesic metric space.
\item[b)]
For each $i\in I$ the intersection
$X_i\cap X_{i+1}$ is a strictly convex closed subset of
both $X_i,X_{i+1}$, and $X_i\cap X_j=\emptyset$
for $\vert i-j\vert \geq 2$.
\item[c)] For each $i\in I$, the sets
$X_i\cap X_{i-1}$ and $X_i\cap X_{i+1}$
are $R$-separated in $X_i$.
\end{enumerate}
Then $X$ is $\delta^\prime$-hyperbolic for a constant
$\delta^\prime>0$ only depending on $\delta,R$.}

{\it Proof:} Let $X,I,X_i$ be as in the lemma.
We may assume without loss of generality that
$I=\mathbb{Z}$. Write $E_i=X_i\cap X_{i+1}$; by
our assumption, $E_i$ is a strictly convex
subset of the proper
$\delta$-hyperbolic spaces $X_i,X_{i+1}$; moreover,
the subsets $E_{i-1},E_i$ of $X_i$ are $R$-separated
for a constant $R>0$ not depending on $i$.
Thus
after possibly enlarging $R$ the following
properties are satisfied.
\begin{enumerate}
\item[i)] Every point $x\in X_i$
can be connected to $E_{i}$
by a geodesic $\zeta_x^+:[0,1]\to X_i$
of minimal length and to $E_{i-1}$ by a geodesic
$\zeta_x^-:[0,1]\to X_i$ of minimal length.
If the distance between $x,y$ is at most $1$ then
the Hausdorff distance between $\zeta_x^{\pm}$ and
$\zeta_y^{\pm}$ is at
most $R$.
\item[ii)] Let $x,y\in X_i$ and let $\gamma$ be a
geodesic connecting $x$ to $y$. Then $\gamma$ is contained in the
$R$-neighborhood of the piecewise geodesic $\tilde \gamma^+$
which is composed of
the arc $\zeta_x^+$, a geodesic in
$E_i$ connecting $\zeta_x^+(1)$ to $\zeta_y^+(1)$ and
the inverse of $\zeta_y^+$. The geodesic $\gamma$ is also
contained in the $R$-neighborhood of
a piecewise geodesic
$\tilde \gamma^-$ which is constructed in the same way using the
geodesic arcs
$\zeta_x^-,\zeta_y^-$ and a geodesic in $E_{i-1}$.
\item[iii)] Let $\gamma_i:[0,1]\to X_i$
be a minimal geodesic connecting $E_{i-1}$
to $E_i$; then for every $x\in
E_{i-1}$ the Hausdorff distance between a minimal geodesic
connecting $x$ to $E_{i}$ and the composition
with $\gamma_i$ of a
geodesic in $E_{i-1}$ connecting $x$ to
$\gamma_i(0)$ is not bigger than $R$.
\end{enumerate}

We use once more the criterion for hyperbolicity
from Lemma 3.3. Namely, we define
in three steps for any pair of
points $x,y\in X$ a curve $c(x,y)$ connecting $x$ to $y$ as
follows.

{\sl Step 1:} If there is some $i\in \mathbb{Z}$ such that $x,y\in
X_i$ then define $c(x,y)$ to be a geodesic in
$X_i$ connecting $x$ to $y$.

{\sl Step 2:} If there is some $i\in \mathbb{Z}$ such that $x\in
X_i-E_i,y\in X_{i+1}-E_i$ then define
$c(x,y)$ to be the piecewise geodesic which is composed from
the geodesic
$\zeta_x^+$
connecting $x$ to $E_i$, a
geodesic arc in $E_i$ connecting $\zeta^+_x(1)$
to $\zeta^-_y(1)$ and the inverse of the geodesic
$\zeta_y^-$.

{\sl Step 3:} If $x\in X_i-E_i,y\in X_j-E_{j-1}$
for some $j\geq i+1$ then define
inductively $c(x,y)$ to be
the piecewise geodesic which consists of the geodesic segment
$\zeta_x^+$, a geodesic in $E_i$ connecting
$\zeta_x^+(1)$ to $\gamma_{i+1}(0)$ and the arc
$c(\gamma_{i+1}(0),y)$.

Assume that the curves $c(x,y)$ are all parametrized on the
unit interval $[0,1]$.
We claim that there is a number $D>0$ only depending on
$\delta,\kappa,n,A$ such that
the curves $c(x,y)$ satisfy the three conditions in
Lemma 3.3.

The first property is immediate from the definition
of the curves $c(x,y)$.
To show that the second condition is valid as well,
let $x,y\in X$, let $0\leq s\leq t\leq 1$ and
let $x^\prime=c(x,y)(s),y^\prime=c(x,y)(t)$ be points on the
curve $c(x,y)$.
We have to show that the Hausdorff distance
between $c(x,y)[s,t]$ and $c(x^\prime,y^\prime)[0,1]$
is bounded from above by
a constant $D_1>0$ only depending on $\delta,\kappa,n,A$.
For this we distinguish three cases.

First, if $x^\prime,y^\prime\in X_i$ for some
$i\in \mathbb{Z}$, then either $c(x,y)[s,t]$ is a geodesic
in $X_i$ connecting
$x^\prime$ to $y^\prime$ or $x^\prime\in E_{i-1}$
or $y^\prime\in E_i$
and by properties ii) and iii) above
for $X_i$,
the Hausdorff distance
between the arc $c(x,y)[s,t]$ and the geodesic
$c(x^\prime,y^\prime)$ connecting
$x^\prime$ to $y^\prime$ is bounded from above by a universal
constant $\chi_1>0$.

Next assume that $x^\prime\in X_i-E_i$ for some
$i\in \mathbb{Z}$ and that $y^\prime\in X_{i+1}-E_i$.
Let $\zeta_{x^\prime}^+$ be a geodesic
of minimal length connecting
$x^\prime$ to $E_i$. By the definition of the
curve $c(x,y)$ and hyperbolicity of
$X_i$,
there is a number $s^\prime>s$ such that
$c(x,y)(s^\prime)\in E_i$ and that
the Hausdorff distance
between $c(x,y)[s,s^\prime]$ and the geodesic $\zeta_{x^\prime}^+$
is at most $R$.
Similarly, by property iii) above and the definition of the
curves $c(v,w)$
there is a number $t^\prime\leq t$ such that
$c(x,y)(t^\prime)\in E_i$ and that
the Hausdorff distance between $c(x,y)[t^\prime,t]$ and
the geodesic $\zeta_{y^\prime}^-$ of minimal length connecting
$y^\prime$ to $E_i$ is bounded from above
by $R$. On the other hand, $c(x^\prime,y^\prime)$
is composed of the arcs $\zeta^+_{x^\prime},\zeta^-_{y^\prime}$
and a geodesic arc in $E_i$ connecting
$\zeta_{x^\prime}(1)$ to $\zeta_{y^\prime}(1)$; moreover,
$c(x,y)[s^\prime,t^\prime]$ is a geodesic in $E_i$
connecting $c(x,y)(s^\prime)$ to $c(x,y)(t^\prime)$.
Since $E_i$ is $\delta$-hyperbolic
for a number $\delta >0$ not depending on $i$, the Hausdorff
distance between any two compact geodesic arcs in $E_i$
is up to an additive constant
bounded from above by the sum of the distances
between the endpoints of the arcs. Therefore
the Hausdorff distance
between $c(x,y)[s,t]$ and $c(x^\prime,y^\prime)$ is at most
$\chi_2$ for a constant $\chi_2\geq \chi_1$ only depending
on $\delta$.

Finally, the case that $x^\prime\in X_i-E_i$ and
$y^\prime\in X_j-E_{j-1}$ for some
$j\geq i+1$ follows immediately from the above consideration.
Namely, in this case there are numbers $s\leq s^\prime<t^\prime \leq t$,
$0\leq \sigma<\tau\leq 1$
such that the arcs $c(x,y)[s^\prime,t^\prime],
c(x^\prime,y^\prime)[\sigma,\tau]$ coincide and that moreover
the above consideration can be applied to the curves
$c(x^\prime,y^\prime)[0,\sigma],
c(x^\prime,y^\prime)[\tau,1]$ and
$c(x,y)[0,s^\prime],c(x,y)[t^\prime,1]$.
Thus the second condition in the proof of Lemma 3.3 is satisfied
for our system of curves with a number $D_1>0$ only
depending on $\delta,\kappa,n,A$ (note that we can choose
$D_1=2\chi_2$ where $\chi_2>0$ is as above).

We are left with showing the thin triangle condition for our
system of curves $c(x,y)$, i.e. we have to find a number $D_2>0$
such that for every triple of points $x,y,z\in Y$ the curve
$c(x,y)$ is contained in the $D_2$-neighborhood of $c(y,z)\cup
c(z,x)$. Consider first the case that the points $x,y,z$ are all
contained in $X_i$ for some $i\in \mathbb{Z}$. Then the curves
$c(x,y),c(y,z),c(z,x)$ are geodesics in $X_i$ connecting these
three points and hence the curve $c(x,y)$ is contained in a
uniformly bounded neighborhood of $c(y,z)\cup c(z,x)$ by
hyperbolicity of $X_i$. Next assume that two of the points, say
the points $x,y$, are contained in $X_i$ but that the third point
$z$ is contained in $X_j-E_{j-1}$ for some $j\geq i+1$. Then the
intersections with $\cup_{p\geq i+1}X_p-E_{i}$ of the curves
$c(x,z),c(y,z)$ coincide. Let $t_x,t_y\in [0,1]$ be such that
$c(x,z)(t_x,1]= c(x,z)[0,1]\cap (\cup_{p\geq i+1}X_p-E_{i})$ and
similarly for $c(y,z)$; then $v=c(x,z)(t_x)=c(y,z)(t_y)$. Together
with property 2 for our curve system established above we conclude
that it is enough to establish the $D_2$-thin triangle condition
for the curves $c(x,v),c(y,v),c(x,y)$. However, since $x,y,v\in
X_i$ this condition holds by our above consideration. The same
argument can also be applied in the case that for each $i$, the
set $X\times [t_{i-1},t_i]$ contains at most one of the points
$x,y,z$. From this we immediately deduce that the third condition
for our curve system is valid as well for a universal constant
$D_2>0$ only depending on $\delta,R$. As a consequence of Lemma
3.3, the space $X$ is $\delta^\prime$-hyperbolic for a constant
$\delta^\prime$ only depending on $\delta,R$. \qed

\bigskip

Let $T$ be a simplicial
tree of bounded valence. Then for any two
points in $T$, there is a unique simple path connecting these
points. For every metric fibration $Y=X\times T\to X$ and
every point $\tau\in T$, the set $X\times \{\tau\}\subset Y$ is
strictly convex. We use these facts together with the
glueing lemma
to extend Lemma 3.4 as
follows.

\bigskip

{\bf Corollary 3.6:} {\it Let $X\times T\to X$ be a bounded metric
fibration with fibre a simplicial tree of bounded valence. Assume
that $X$ is $\delta$-hyperbolic for some $\delta >0$ and that
vertical distances satisfy the $\kappa,n$-flaring property with
threshold $A>0$ for some $\kappa>1,n>0,A>0$. Then $Y$ is
$\delta_1$-hyperbolic for a number $\delta_1>0$ only depending on
$\kappa,n,\delta,A$.}

{\it Proof:} We begin with showing the corollary in the particular
case that the tree $T$ is just an arbitrary closed connected
subset $J$ of
the real line $\mathbb{R}$. Thus let $X$ be a $\delta$-hyperbolic
geodesic metric space, let $J\subset \mathbb{R}$ be an arbitrary
closed connected set and let $Y=X\times J$ be a metric fibration
with fibre $J$. Assume that vertical distances satisfy the
$\kappa,n$-flaring property with threshold $A$ for some $\kappa
>1, n>0,A>0$ and assume without loss of generality that $A\geq 1$.

Let $0\in J$ and assume that $0$ is an endpoint of $J$ if
$J\not=\mathbb{R}$. Assume moreover that in this case the set $J$
is contained in $[0,\infty)$. For $t\in J$ let $\ell_t^1:X\to
[0,\infty)$ be the function which associates to a point $x\in X$
the length of the vertical path $\{x\}\times [0,t]$. By
assumption, the function $\ell^1_t$ is continuous. Write
$\mu^1(t)=\inf_{x\in X} \ell^1_t(x)$; the function $t\to \mu^1(t)$
is continuous and monotonously increasing on $[0,\infty)$,
monotonously decreasing on $(-\infty,0]$. Let $t_1\in (0,\infty]$
be the smallest positive number with $\mu^1(t_1)=A$; here we write
$t_1=\infty$ if $\mu^1(t)<A$ for all $t>0$. If $t_1<\infty$ then
define for $t\geq t_1$ a new function $\ell^2_t: X\to [0,\infty)$
by assigning to $x\in X$ the length of the arc $\{x\}\times
[t_1,t]$. Let $\mu_2(t)=\inf_{x\in X}\ell_t^2(x)$ and let $t_2\in
(t_1,\infty]$ be the smallest number such that $\mu^2(t_2)=A$.
Inductively we construct in this way an increasing sequence
$0<t_1<t_2<\dots$ and functions $\mu^i,\ell_t^i$. The sequence
might be trivial, finite or infinite. If $J=\mathbb{R}$ then
define in the same way a sequence $0>t_{-1}>t_{-2}>\dots$ and
functions $\mu^i,\ell_t^i$ $(i\leq -1)$.

By Lemma 3.4, there is a constant $\delta_0>0$ such that for each
$i\in \mathbb{Z}$, the convex subset $X\times [t_{i-1},t_i]$ of
$X\times J$ is $\delta_0$-hyperbolic. The sets $X\times
\{t_{i-1}\},X\times \{t_i\}$ are strictly convex in $Y$. Moreover,
by the remark preceding Lemma 3.5 they are also $R$-separated for
a constant $R>0$ only depending on $\kappa,n,A$. Thus we can apply
Lemma 3.5 and conclude that the metric fibration $X\times J$ is
$\delta_1$-hyperbolic for a constant $\delta_1>0$ only depending
on $\delta,\kappa,n,A$.

Now let $T$ be a simplicial tree of bounded valence.
Let $X\times T\to X$ be a metric fibration over a proper
$\delta$-hyperbolic geodesic metric space $X$. Assume that
vertical distances satisfy the $\kappa,n$-flaring property
with threshold $A>0$ for some $\kappa>1,n>0,A>0$.
Our goal is to show that $Y$ is $\delta_2$-hyperbolic
for a constant $\delta_2>0$ only depending on
$\delta,\kappa,n,A$.

For this let $y_1,y_2,y_3\in Y$ be a triple of points and
let $c(y_i,y_j)$ $(i,j=1,2,3)$ be geodesics in $X\times T$
connecting $y_i$ to $y_j$. We have to show that
$c(y_1,y_2)$ is contained in a uniformly bounded neighborhood of
$c(y_2,y_3)\cup c(y_3,y_1)$.
Write $y_i=(x_i,\tau_i)$ with $x_i\in X,\tau_i\in T$.
For $i=1,2,3$ let $J_i$ be the unique embedded segment in
$T$ connecting $\tau_i$ to $\tau_{i+1}$ (indices are taken mod 3).
Then the intersection $\cap_i J_i$ consists of a unique point
$\tau$. By our above consideration, the subsets
$X\times J_i\subset Y$ of $Y$ are strictly convex and
moreover $\delta_1$-hyperbolic for a universal constant
$\delta_1>0$; they contain $X\times \{\tau\}$
as a strictly convex subset.
Let $\rho(y_i,y_{i+1})$ be
a piecewise geodesic which up to orientation and parametrization
is composed of a minimal geodesic
$\alpha_i:[0,1]\to X\times J_i$
connecting $y_i$ to $X\times \{\tau\}$, a minimal geodesic
$\alpha_{i+1}:[0,1]\to X\times J_{i+1}$ connecting
$y_{i+1}$ to $X\times \{\tau\}$ and a geodesic in $X\times \{\tau\}$
connecting $\alpha_i(1)$ to $\alpha_{i+1}(1)$.
The Hausdorff distance between
the geodesic $c(y_i,y_{i+1})\subset
X\times J_i$ and the piecewise geodesic $\rho(y_i,y_{i+1})$
is bounded
from above by a constant only depending on $\delta,\kappa,n,A$.
Since $X\times \{\tau\}$ is $\delta$-hyperbolic, from this
hyperbolicity of $Y$ is immediate.
This shows the corollary. \qed

\bigskip

We summarize the results of this section as follows.

\bigskip

{\bf Corollary 3.7:} {\it Let $X$ be a proper hyperbolic geodesic
metric space and let $Y=X\times T\to X$ be a bounded metric
fibration with fibre a simplicial tree of bounded valence. Then
$Y$ is hyperbolic if and only if vertical distances satisfy the
$\kappa,n$-flaring property with threshold $A$ for some $\kappa
>1,n>0,A>0$.}

{\it Proof:} Let $X$ be a proper hyperbolic geodesic metric space
and let $Y=X\times T\to X$ be a bounded metric fibration with
fibre a simplicial tree of bounded valence. Lemma 5.2 of \cite{FM}
shows that if $Y$ is hyperbolic, then vertical distances satisfy
the $\kappa,n$-flaring property with threshold $A$ for some
$\kappa>0,n>0,A>0$. By Corollary 3.6, this condition is
also sufficient for hyperbolicity of $Y$. \qed

\section{Proof of the theorem}

In this final section we consider
a \emph{closed} surface
of genus $g\geq 2$.
Our goal is to show that for a convex cocompact subgoup $\Gamma$
of the mapping class group
${\cal M}_g$ for $S$, the natural $\pi_1(S)$-extension $\Gamma_S$ of
$\Gamma$ is word hyperbolic. For this choose a finite symmetric
generating set ${\cal G}$ for $\Gamma$ and denote by $\Vert
\,\Vert$ the induced word norm on $\Gamma$ and by ${\cal C\cal G}$
the corresponding Cayley graph. Choose a point $h$
in the Teichm\"uller space ${\cal T}_g$ for $S$
which does not admit any nontrivial automorphisms (recall that the
set of such points is open and dense in ${\cal T}_g$) and define a
map $\Theta:{\cal C\cal G}\to {\cal T}_g$ by mapping a vertex
$\phi\in \Gamma$ to the point $\phi h\in {\cal T}_g$ and by
mapping an edge $e$ of ${\cal C\cal G}$ to the Teichm\"uller
geodesic arc connecting the image of the endpoints of $e$. By
Lemma 2.8, the map $\Theta$ is 
a quasi-isometric embedding; moreover, the set
$\Theta {\cal C\cal G}$ is invariant under the action of $\Gamma$.

There is a natural smooth marked
surface bundle ${\cal
S}\to {\cal T}_{g}$ whose fibre ${\cal S}_z$
at a point $z\in {\cal T}_g$ is just the surface
$S$ with the marking defined by $z$. The
hyperbolic structure $z\in {\cal T}_g$
defines a smooth Riemannian
metric on the fibre ${\cal S}_z$, and these metrics fit
together to a smooth Riemannian metric on the
\emph{vertical bundle} of
the fibration (i.e. the tangent bundle of the fibres).
In other words, the vertical foliation of ${\cal S}$ into
the fibres of our fibration admits a natural smooth
Riemannian metric.

The action of the mapping class group ${\cal M}_g$ on ${\cal T}_g$
lifts to a unique action
on ${\cal S}$ which is determined by the
requirement that for every $\phi\in
{\cal M}_g$ and every $z\in {\cal T}_g$,
the restriction of the lift of $\phi$ to ${\cal S}_z$
is the unique isometry of ${\cal S}_z$ onto
${\cal S}_{\phi z}$ in the isotopy class
determined by $\phi$. In particular, the Riemannian metric on the
vertical foliation is invariant under the action of ${\cal M}_g$.
The restriction ${\cal S}_\Gamma$ of the bundle ${\cal S}$ to
$\Theta{\cal C\cal G}$ is invariant under the action of
the subgroup $\Gamma$ of ${\cal M}_g$.

We equip now the bundle ${\cal S}_\Gamma$ with the following
geodesic metric. First, recall that for a given edge $b$ in
${\cal C\cal G}$ the arc $\Theta b$ is a geodesic, and its
endpoints are marked hyperbolic metrics on the surface $S$ which
are isometric with an isometry in the class determined by the
element of ${\cal G}$ corresponding to $b$. If we identify the
edge $b$ with the unit interval $[0,1]$ then for all $s,t\in
[0,1]$ there is a unique Teichm\"uller map of minimal
quasi-conformal dilatation which maps the fibre ${\cal
S}_{\Theta(s)}$ to ${\cal S}_{\Theta(t)}$, and these maps combine
to a smooth fibre preserving horizontal flow on the restriction
${\cal S}_{\Theta b}$ of ${\cal S}$ to $\Theta b$. Defining the
tangent of each of these flow-lines to be orthogonal to the fibres
and of the same length as its projection to $\Theta b$ defines a
smooth Riemannian metric on ${\cal S}_{\Theta b}$ so that the
canonical projection ${\cal S}_{\Theta b}\to \Theta b$ is a
Riemannian submersion. If two edges are incident on the same
vertex, then the metrics on the fibres over this vertex coincide.
Therefore, the metrics naturally induce a
complete length metric $d$ on
${\cal S}_\Gamma$ which for every edge $b$ of ${\cal C\cal
G}$ restricts to the length metric of the above Riemannian
structure on ${\cal S}_{\Theta b}$.

The universal cover ${\cal H}$ of ${\cal S}$ is a smooth fibre
bundle $\Pi:{\cal H}\to {\cal T}_g$ whose fibre ${\cal H}_z$ at a
point $z\in {\cal T}_g$ equipped with the lift of the Riemannian
metric on ${\cal S}_z$ is isometric to the hyperbolic plane ${\bf
H}^2$. The group $\Gamma_S$ acts on ${\cal H}$ as a group of
bundle isomorphisms preserving the metric of the vertical
foliation. The pre-image of every $\Gamma$-invariant subset of
${\cal T}_g$ is invariant under the action of $\Gamma_S$. In
particular, the set ${\cal H}_\Gamma=\Pi^{-1}(\Theta {\cal C\cal
G})$ is $\Gamma_S$-invariant. The metric on ${\cal S}_\Gamma$
lifts to a geodesic metric $d$ on ${\cal H}_\Gamma$. The group
$\Gamma_S$ acts on the geodesic metric space $({\cal H}_\Gamma,d)$
isometrically, properly and cocompactly and hence we have (see
\cite{FM}).

\bigskip

{\bf Lemma 4.1:} {\it ${\cal H}_\Gamma$ is quasi-isometric to
$\Gamma_S$.}

\bigskip

By Lemma 4.1 it is therefore enough to show that
the bundle ${\cal H}_\Gamma$ with its $\Gamma_S$-invariant
geodesic metric
is hyperbolic.

To show that this is indeed the case, we use the
results from Section 3, applied to suitably defined
line-subbundles of ${\cal H}_\Gamma$. For the construction
of these bundles,
fix a standard system $a_1,b_1,\dots,a_g,b_g$ of generators for
the fundamental group $\pi_1(S)$ of $S$. For every $z\in
{\cal T}_g$ there is a unique isomorphism $\rho(z)$ of
$\pi_1(S)$ onto a discrete subgroup $\Upsilon(z)$ of
$PSL(2,\mathbb{R})$ such that the surface ${\bf H}^2/\Upsilon(z)$
is isometric to ${\cal S}_{z}$ and that moreover the following
holds (see \cite{IT}).
\begin{enumerate}
\item[a)] The conjugacy class of the representation
$\rho(z)$ is determined by the marking of ${\cal S}_z$.
\item[b)]
In the upper half-plane model for ${\bf H}^2$, the points
$0,\infty$ are attracting and repelling fixed points for the
action of $\rho(z)(b_g)$, and the point $1$ is an attracting fixed
point for the action of $\rho(z)(a_g)$.
\end{enumerate}
The representation $\rho(z)$ depends smoothly on $z$, and
the surface bundle ${\cal S}$ over ${\cal T}_g$ is the quotient
of the trivial bundle ${\cal T}_g\times {\bf H}^2$ under
the action of the group $\pi_1(S)$ defined by
$\phi(z,v)=(z,\rho(z)(\phi)(v))$ $(\phi\in \pi_1(S),(z,v)\in
{\cal T}_g\times {\bf H}^2$).

For every $z\in
{\cal T}_g$ the fibre ${\cal H}_z$ of the bundle ${\cal H}$
admits a compactification by adding the \emph{ideal
boundary} $\partial {\cal H}_z$. Every pair of distinct points in
$\partial{\cal H}_z$ defines uniquely a geodesic line in ${\cal
H}_z$. Let again $h\in {\cal T}_g$ be a point whose
$\Gamma$-orbit is the vertex set of $\Theta{\cal C\cal G}$.
For every $z\in \Theta{\cal C\cal G}$, the isomorphism
$\rho(z)\circ \rho(h)^{-1}$ of $\Upsilon(h)$ onto
$\Upsilon(z)$
induces a homeomorphism $\omega(z)$ of
$\partial{\cal H}_h$ onto $\partial{\cal H}_z$. For every
pair of distinct points
$\xi\not=\eta\in \partial{\cal H}_h$ we define
a line subbundle
${\cal L}^{\xi,\eta}$ of ${\cal H}_\Gamma$ by requiring that
its fibre ${\cal L}^{\xi,\eta}_z$
at $z$ is the geodesic line in ${\cal H}_z$ whose endpoints
in $\partial{\cal H}_z$ are the images of the points $\xi,\eta$
under the homeomorphism $\omega(z)$.
We equip ${\cal L}^{\xi,\eta}$ with a complete length metric
whose restriction to each fibre coincides with the
restriction of the metric on ${\cal H}$ and which
is such that the following holds. For each $z\in \Theta{\cal C\cal G}$
let ${\cal B}_z\subset{\cal H}_z$ be the tubular neighborhood of
radius one about the fibre ${\cal L}^{\xi,\eta}_z$ of
${\cal L}^{\xi,\eta}$ at $z$; then
${\cal B}=\cup_{z\in \Gamma} {\cal B}_z$ is
a fibre bundle over $\Theta{\cal C\cal L}$ which is an open
subset of ${\cal H}_\Gamma$. The restriction to ${\cal B}$
of the length structure on ${\cal H}_\Gamma$ induces
a length metric on ${\cal B}$; we require that the
inclusion ${\cal L}^{\xi,\eta}\to {\cal B}$ is a
$p$-quasi-isometry for some constant $p\geq 1$ 
independent of $\xi,\eta$.
We have.

\bigskip

{\bf Lemma 4.2:} {\it There is a number $q>0$ so that for every
pair $\xi\not=\eta\in \partial {\cal H}_h$ the following is
satisfied.
\begin{enumerate}
\item The inclusion ${\cal L}^{\xi,\eta}\to {\cal H}_\Gamma$ is a
$q$-quasi-isometric embedding.
\item The bundle ${\cal L}^{\xi,\eta}$ is $q$-hyperbolic.
\end{enumerate}
}

{\it Proof:} Let $\xi\not=\eta\in \partial{\cal H}_h$. We
begin with showing
that there is a number $q_0>0$ not depending on $\xi,\eta$ so that
the inclusion $\iota:{\cal L}^{\xi,\eta}\to {\cal H}_\Gamma$ is a
$q_0$-quasi-isometric embedding. For this we have to show that for
any two points $v,w\in {\cal L}^{\xi,\eta}$ the distance
in ${\cal L}^{\xi,\eta}$ between
$v,w$ is not bigger than $q_0$ times the distance between
$\iota(v),\iota(w)$. By definition, for the neighborhood
${\cal B}$ of $\iota({\cal L}^{\xi,\eta})$ in ${\cal H}_\Gamma$
which intersects each fibre ${\cal H}_z$ in the
neighborhood of radius one about the geodesic
$\iota{\cal L}_z^{\xi,\eta}$, we have to find
a curve connecting $\iota(v)$ to $\iota(w)$ in ${\cal B}$ 
whose length
is bounded from above by a constant multiple of the
distance between $\iota(v)$ and $\iota(w)$ in ${\cal H}_\Gamma$.

For this note first that by construction of the
metric on ${\cal H}_\Gamma$, there is
a universal number $a>0$ with the following property.
If $y\in \iota({\cal L}^{\xi,\eta})$ and if $\zeta:[0,a]\to
{\cal H}_\Gamma$ is a horizontal curve of length at most
$a$ issuing from $y$, then $\zeta[0,a]\subset {\cal B}$.

Let $P:{\cal H}_\Gamma \to{\cal L}^{\xi,\eta}$ be the
unique bundle map whose restriction to a fibre ${\cal H}_z$ is the
shortest distance projection of ${\cal H}_z$ onto ${\cal
L}^{\xi,\eta}_z$. Let $\zeta:[0,am]\to {\cal H}_\Gamma$ be any
geodesic of length $am\geq 0$
connecting the points $\iota(v),\iota(w)\in {\cal
L}^{\xi,\eta}$. Since there is a number $n>0$ such that
the horizontal transport of
the fibres of the bundle ${\cal H}\to \Theta{\cal C\cal G}$
along a geodesic arc in $\Theta{\cal C\cal G}$
of length at most $\ell$
is a bilipschitz map with bilipschitz
constant bounded from above by $e^{n\ell}$,
there is a curve $\zeta_0:[0,2m]\to
{\cal H}_\Gamma$ connecting $\iota(v)$ to $\iota(w)$
with the
property that for every $i< m$ the restriction of
$\zeta_0$ to the interval $[2i,2i+1]$ is
a horizontal arc of length at most $a$,
and the restriction of $\zeta_0$ to $[2i+1,2i+2]$ is
vertical and of length at most $e^{an}$.
The length of $\zeta_0$ is bounded from above
by $am(e^{an}+1)$.
Define a curve
$\zeta_1:[0,2m]\to {\cal H}_\Gamma$ by requiring that for each
$i<m$, the restriction of $\zeta_1$ to $[2i,2i+1]$ is the
horizontal arc issuing from $P\zeta_0(2i)$ whose projection to
$\Theta{\cal C\cal G}$ coincides with the projection of
$\zeta_0[2i,2i+1]$, and the restriction of $\zeta_1$ to
$[2i+1,2i+2]$ is the vertical geodesic which connects the point
$\zeta_1(2i+1)$ to $P\zeta_0(2i+2)$. Note that
this curve is entirely contained in ${\cal B}$.
To establish our claim it is now enough to show that the
distance between the points $\zeta_1(2i+1)$ and
$P\zeta_0(2i+2)$ is bounded from above by a universal constant.
Namely, if this is the case then the length
of $\zeta_1$ is bounded from above by a universal
multiple of the length of the geodesic $\zeta$.

For this fix for the moment an arbitrary number $L>1$.
Let $\Psi:{\bf H}^2\to {\bf H}^2$ be an
$L$-quasi-isometry which induces the homeomorphism $\psi$
of the ideal
boundary $\partial {\bf H}^2$ of ${\bf H}^2$.
Let
$\gamma:\mathbb{R}\to {\bf H}^2$ be a geodesic line and
let $R:{\bf H}^2\to \gamma$ be the shortest distance projection.
Let $y\in {\bf H}^2$ be such that $R(y)=\gamma(0)=x$ and let
$\zeta:[0,\infty)\to {\bf H}^2$ be the geodesic ray issuing
from $\zeta(0)=x$ and passing through $y$. Then the
concatenation $\sigma_{\pm}$ of the inverse
$\zeta^{-1}$ of $\zeta$
with the geodesic ray $\gamma[0,\infty)$ and
the inverse of the geodesic ray $\gamma(-\infty,0]$
is a uniform quasi-geodesic in ${\bf H}^2$ containing both
$x$ and $y$. 
Since the isometry group of ${\bf H}^2$ acts triply transitive
on the ideal boundary,
via composing $\Psi$ with an isometry we may assume
that $\psi$ fixes the endpoints $\gamma(\pm \infty),
\zeta(\infty)$ of $\gamma,\zeta$.
Now $\Psi$ is an $L$-quasi-isometry and hence
the image $L\gamma$ of $\gamma$ is an $L$-quasi-geodesic
with the same endpoints as $\gamma$. Thus the Hausdorff
distance between $\gamma$ and $\Psi\gamma$ is uniformly
bounded and hence the
distance between $\Psi x$ and $R\Psi x$ is bounded from above
by a universal constant $p_0>0$. Similarly,
$\Psi\sigma_{\pm}$ are uniform quasi-geodesics
contained in a uniformly bounded neighborhood of $\sigma_{\pm}$.
Since $\Psi x$ is contained in $\Psi\sigma_+\cap \Psi\sigma_-\cap
\Psi\gamma$, we conclude that the distance between
$x$ and $\Psi x$ is uniformly bounded and that the distance
between $R\Psi y$ and $\Psi x$ is uniformly bounded as well.

Now for every arc $\nu:[0,1]\to \Theta{\cal C\cal G}$
of length at most $a$,
the homeomorphism $\zeta:{\cal H}_{\nu(0)}\to {\cal H}_{\nu(1)}$
obtained by horizontal transport of the fibres along $\nu$
is an $L$-quasi-isometry for a universal constant $L>1$
which induces the boundary homeomorphism
$\omega(\nu(1))\circ\omega(\nu(0))^{-1}$.
By the definition of the line bundle ${\cal L}^{\xi,\eta}$ and
the above consideration, we conclude that the distance
in ${\cal H}_{\zeta_0(2i+1)}$ between the points
$\zeta_1(2i+1)$ and the point $P\zeta_0(2i+1)$ is uniformly bounded.
As a consequence, the inclusion $\iota:{\cal L}^{\xi,\eta}\to
{\cal H}_\Gamma$ is a $q_0$-quasi-isometric
embedding for a constant $q_0>0$ not depending on $\xi,\eta$.

Next we claim that ${\cal L}^{\xi,\eta}$ is uniformly
quasi-isometric to a metric fibration over $\Theta{\cal C\cal G}$
with fibre $\mathbb{R}$. Namely, fix a component $A$ of
$\partial{\cal H}_{h}-\{\xi,\eta\}$. For $z\in
\Theta{\cal C\cal G}$ and $\nu\in A$ define
$\Pi(\nu,z)\in {\cal L}^{\xi,\eta}_z$ to be the shortest
distance projection of $\rho(z)(\nu)$ to ${\cal L}^{\xi,\eta}_z$.
The map $\nu\to \Pi(\nu,z)$ is a homeomorphism of $A$ onto ${\cal
L}^{\xi,\eta}_z$. It follows from
our above consideration that for
every $\phi\in {\cal G}$ and every $\nu\in A$ the distance
between $\Pi(h,\nu)\in {\cal H}_h$ and
$\Pi(\phi h,\nu)\in {\cal H}_{\phi h}$ is uniformly
bounded. But the map $\Pi$ is equivariant with
respect to the action of
$\Gamma$ on $\Theta{\cal C\cal G}$ and on ${\cal H}_\Gamma$
and hence for every
$\psi\in \Gamma$, the distance between
$\Pi(\phi\psi h,\nu)$ and $\Pi(\psi h,\nu)$ is
bounded from above by
the same constant. As a consequence,
for every fixed $\nu\in A$ the
map $z\to \Pi(\nu,z)$ is a $q_1$-quasi-isometric embedding
of $\Theta{\cal C\cal L}$ into ${\cal L}^{\xi,\eta}$ for
a number $q_1>0$ not depending on $\nu$ and on $\xi,\eta$,
and these quasi-isometric embeddings can be used to define
on ${\cal L}^{\xi,\eta}$ the structure of a bounded metric
fibration which is quasi-isometric to ${\cal L}^{\xi,\eta}$
equipped with the
metric induced from the metric on ${\cal H}_\Gamma$.

By Corollary 3.6, to show that ${\cal L}^{\xi,\eta}$ is
$q$-hyperbolic for a constant $q>0$ not depending on $\xi,\eta$ we
only have to show that vertical distances for our metric fibration
satisfy the $\kappa,n$-flaring property with threshold $A$ for
numbers $\kappa>1,n>0,A>0$ not depending on $\xi,\eta$.

For this let $\zeta:\mathbb{R}\to \Theta{\cal C\cal G}$ be any
geodesic line. By the discussion in Section 2, $\zeta$ is
contained in a uniformly bounded neighborhood of a Teichm\"uller
geodesic. By the results of Mosher \cite{Mo03},
the restriction ${\cal H}_\zeta$
of the bundle ${\cal H}_\Gamma$ to $\zeta$ is
$\delta_0$-hyperbolic for a constant $\delta_0$ not depending on
$\zeta$. The above consideration can be applied to the
restrictions of the bundles ${\cal L}^{\xi,\eta}$ and ${\cal
H}_\Gamma$ to the geodesic $\zeta$ and shows that the inclusion
${\cal L}^{\xi,\eta}\vert \zeta\to {\cal H}_\zeta$ is a
quasi-isometric embedding. Since ${\cal H}_\zeta$ is hyperbolic,
the bundle ${\cal L}^{\xi,\eta}\vert \zeta$ is
$\delta_1$-hyperbolic for a universal constant $\delta_1$. Now
the geodesic $\zeta$ was arbitrary and therefore Lemma 5.2 of
\cite{FM} shows that vertical distances in ${\cal L}^{\xi,\eta}$ satisfy
the $\kappa,n$-flaring property with threshold $A$ for constants
$\kappa>1,n>0,A>0$ not depending on $\xi,\eta$. As a consequence of
this and Corollary 3.6, the line bundle ${\cal L}^{\xi,\eta}\to
\Theta{\cal C\cal G}$ is $q$-hyperbolic for a universal constant
$q>0$. \qed

\bigskip

Now we are ready to show.

\bigskip

{\bf Lemma 4.3:} {\it The bundle ${\cal H}_\Gamma$ is hyperbolic.}

{\it Proof:} Fix a triple of pairwise distinct points
$\xi_1,\xi_2,\xi_3\in \partial{\cal H}_h$. Then the line bundles
${\cal L}^{\xi_i,\xi_{i+1}}$ bound a subbundle ${\cal V}$ of
${\cal H}_\Gamma$ whose fibre at a point $z$ is isometric to an
ideal triangle in ${\cal H}_z$. The arguments in the proof of
Lemma 4.2 show that for a suitable choice of a length metric
on ${\cal V}$ which restricts to the metric on the fibres
induced from the metrics on ${\cal H}_\Gamma$
the inclusion ${\cal V}\to {\cal H}_\Gamma$ is
a quasi-isometric embedding. We claim that the bundle ${\cal V}$
is $\delta_0$-hyperbolic for a number $\delta_0$ not depending on
$\xi_i$. Namely, an ideal hyperbolic triangle $T$ is uniformly
quasi-isometric to the tripod which consists of the unique point
in the interior of $T$ of equal distance to each of the three sides
and three geodesic rays issuing from this point
which make a mutual angle $2\pi/3$.  The arguments in the proof of
Lemma 4.2 show that the bundle ${\cal V}$ is 
quasi-isometric to a metric
fibration over $\Theta{\cal C\cal C}$ whose fibre is precisely
this tripod.

Now as before, for every geodesic $\zeta$ in $\Theta{\cal C\cal
G}$ the restriction of ${\cal H}_\Gamma$ to $\zeta$ is hyperbolic
and hence the same is true for the restriction ${\cal V}_\zeta$
of the bundle
${\cal V}$ since ${\cal V}_\zeta\subset{\cal H}_\zeta$
is quasi-isometrically
embedded. As a consequence of this and Lemma 5.2 of \cite{FM},
vertical distances in the tripod bundle satisfy the
$\kappa,n$-flaring property with threshold $A$ for universal
numbers $\kappa>1,n>0,A>0$. Corollary 3.6 then shows that the
bundle ${\cal V}$ is $\delta_0$-hyperbolic for a universal number
$\delta_0>0$.

We now use once more Lemma 3.3. Namely, for
$x,y\in {\cal H}_\Gamma$ construct a curve $c(x,y)$ connecting $x$
to $y$ as follows. Fix once and for all a point $\xi\in \partial
{\cal H}_h$. Then $\cup_{\nu\not=\xi}{\cal
L}^{\xi,\nu}_z={\cal H}_z$ for every $z\in
\Theta{\cal C\cal G}$ and therefore
there are unique not necessarily distinct points
$\nu_x,\nu_y\in \partial {\cal H}_h-\{\xi\}$ such that $x\in {\cal
L}^{\xi,\nu_x},y\in {\cal L}^{\xi,\nu_y}$.
The points $\xi,\nu_x,\nu_y$
define a bundle over $\Theta{\cal C\cal G}$ whose fibre is a
(possibly degenerate) ideal
hyperbolic triangle; this bundle is quasi-isometrically embedded
in ${\cal H}_\Gamma$. We define $c(x,y)$ to be a geodesic
in this bundle
connecting $x$ to $y$.

We claim that the system of curves $c(x,y)$
satisfies the properties listed
in Lemma 3.3. Namely, property 1) is immediate from
the fact that the bundles ${\cal V}\subset{\cal H}_\Gamma$
as above are uniformly
quasi-isometrically embedded. To show property 3) above, let
$x,y,z$ be a triple of points. Then $x,y,z$ are contained in a
subbundle of ${\cal H}_\Gamma$ whose fibre at the point $z$
is the ideal quadrangle in ${\cal H}_z$ with vertices
$\xi,\nu_x,\nu_y,\nu_z$.
As before, this bundle is $\delta_1$-hyperbolic for a universal
number $\delta_1>0$, and the subbundles whose fibres are
the ideal triangles with vertices
$\xi,\nu_x,\nu_y$ and
$\xi,\nu_x,\nu_z$ and $\xi,\nu_y,\nu_z$
are uniformly quasi-isometrically
embedded. By definition of our curve system,
property 3) for our curve system now
follows from hyperbolicity of our bundle of quadrangles, and
property 2) is obtained in the same way. \qed

\bigskip

As a consequence, we obtain the second part of our theorem.

\bigskip

{\bf Corollary 4.4:} {\it For a finitely generated subgroup
$\Gamma<{\cal M}_g$, the following are equivalent.
\begin{enumerate}
\item The natural $\pi_1(S)$-extension of $\Gamma$ is hyperbolic.
\item $\Gamma$ is convex cocompact.
\end{enumerate}
}

{\it Proof:} Farb and Mosher \cite{FM} show the following.
If $\Gamma<{\cal M}_g$ is any finitely generated group such that
the $\pi_1(S)$-extension of $\Gamma$ is hyperbolic,
then there is a quasi-convex orbit for the action of
$\Gamma$ on ${\cal T}_g$. By Theorem 2.9, this is equivalent to
saying that $\Gamma$ is convex cocompact.

On the other hand, by Lemma 4.1 and Lemma 4.3
the $\pi_1(S)$-extension
of a convex cocompact
subgroup $\Gamma$ of ${\cal M}_g$ is word hyperbolic.
\qed

\bigskip

{\bf Remark:} As mentioned earlier, there is no
example known of a convex cocompact subgroup of
${\cal M}_g$ which is not virtually free.
On might ask whether indeed \emph{all} convex
cocompact subgroups $\Gamma$ of ${\cal M}_g$ are virtually
free. One possible approach to study this question
is via the Gromov boundary $\partial \Gamma$
of $\Gamma$. Namely, by Theorem 2.9, the Gromov
boundary of $\Gamma$ embeds into the Gromov
boundary of the complex of curves ${\cal C}(S)$, and
its image is contained in the subset of all
minimal geodesic laminations which are \emph{uniquely
ergodic}, i.e. which support up to multiple a
unique transverse invariant measure.

As a consequence, $\partial \Gamma$ embeds into the set ${\cal
U\cal E}\subset {\cal P\cal M\cal L}$ of uniquely ergodic
projective measured laminations. Therefore, if ${\cal U\cal E}$ is
totally disconnected, then the same is true for the Gromov
boundary of $\Gamma$ and consequently $\Gamma$ is virtually free.
However, to my knowledge, nothing is known about the structure of
the set ${\cal U\cal E}$.

\bigskip

{\bf Acknowledgement:} I am grateful to Benson Farb for
inspiration and to Maciej Zworski for hospitality.


\bigskip




\begin{thebibliography}{GDH99}






\bibitem[BF92]{BF} M.~Bestvina, M.~Feighn, {\em
A combination theorem for negatively curved groups},
J. Diff. Geom. 35 (1992), 85--101.





\bibitem[B02]{B1} B.~Bowditch, {\em Intersection numbers
and the hyperbolicity of the curve complex}, preprint 2002.





\bibitem[BH99]{BH} M.~Bridson, A.~Haefliger, {\sl Metric
spaces of non-positive curvature}, Springer Grund\-leh\-ren 319,
Springer, Berlin 1999.




\bibitem[Bu92]{Bu} P.~Buser, {\sl Geometry and spectra
of compact Riemann surfaces}, Birkh\"auser, Boston 1992.



\bibitem[CEG87]{CEG} R.~Canary, D.~Epstein, P.~Green,
{\em Notes on notes of Thurston}, in ``Analytical and geometric
aspects of hyperbolic space'', edited by D.~Epstein, London Math.
Soc. Lecture Notes 111, Cambridge University Press, Cambridge
1987.





\bibitem[FM02]{FM} B.~Farb, L.~Mosher, {\em Convex
cocompact subgroups of mapping class groups},
Geom. Top. 6 (2002), 91-152.



\bibitem[FLP91]{FLP}
A.~Fathi, F.~Laudenbach, V.~Po\'enaru, {\sl Travaux de
Thurston sur les surfaces,} Ast\'erisque 1991.










\bibitem[GDH99]{GDH99} G.~Gonzalez-D\'iez, W.~Harvey,
{\em Surface subgroups inside mapping class groups},
Topology 38 (1999), 57--69.






\bibitem[H04]{H04} U.~Hamenst\"adt, {\em Train tracks and
the Gromov boundary of the complex of curves}, to appear
in ``Spaces of Kleinian groups'', edited by
Y.~Minski, M.~Sakuma and C.~Series, Lond.
Math. Soc. Lec. Notes 329, 187--207;
arXiv:math.GT/0409611.






\bibitem[H05]{H05} U.~Hamenst\"adt, {\em
Geometry of the complex of curves and of Teichm\"uller
space}, arXiv:math.GT/0502256.



\bibitem[Ha81]{Ha} W.~J.~Harvey, {\em Boundary structure
of the modular group}, in ``Riemann Surfaces and Related
topics: Proceedings of the 1978 Stony Brook Conference''
edited by I.~Kra and B.~Maskit, Ann. Math. Stud. 97,
Princeton, 1981.




\bibitem[IT99]{IT} Y.~Imayoshi, M.~Taniguchi,
{\sl An introduction to Teichm\"uller spaces},
Springer, Tokyo 1999.



\bibitem[I02]{I} N.~V.~Ivanov, {\em Mapping class groups},
Chapter 12 in ``Handbook of Geometric Topology'', edited by
R.J.~Daverman and R.B.~Sher, Elsevier Science (2002), 523--633.



\bibitem[KL05]{KL05} R.~Kent, C.~Leininger,
{\em Shadows of mapping class groups: capturing
convex cocompactness}, arXiv:math.GT/0505114.


\bibitem[Ke92]{Ke} S.~Kerckhoff, {\em Lines
of minima in Teichm\"uller space}, Duke Math. J. 65 (1992),
187-213.



\bibitem[K99]{K} E.~Klarreich, {\em The boundary at infinity of the curve
complex and the relative Teichm\"uller space,} preprint 1999.

\bibitem[LR05]{LR05} C.~Leininger, A.~Reid, {\em A combination
theorem for Veech subgroups of the mapping class group},
arXiv:math.GT/0410041.








\bibitem[Ma82a]{Ma82a} H.~Masur, {\em Interval exchange
transformations and measured foliations},
Ann. Math. 115 (1982), 169-200.



\bibitem[Ma82b]{Ma82b} H.~Masur, {\em Two boundaries of
Teichm\"uller space}, Duke Math. J. 49 (1982), 183--190.








\bibitem[MM99]{MM1} H.~Masur, Y.~Minsky, {\em Geometry of the
complex of curves I: Hyperbolicity}, Invent. Math. 138 (1999),
103--149.







\bibitem[MP89]{MP89} J.~McCarthy, A.~Papadopoulos,
{\em Dynamics on Thurston's sphere of projective
measured foliations}, Comm. Math. Helv. 64 (1989), 133--166.









\bibitem[Mo96]{Mo96} L.~Mosher, {\em Hyperbolic
extensions of groups}, J. Pure Appl. Alg.110 (1996),
305--314.



\bibitem[Mo03]{Mo03} L.~Mosher, {\em Stable Teichm\"uller
quasigeodesics and ending laminations}, Geom. Top. 7
(2003), 33-90.


























\end{thebibliography}
\end{document}